\documentclass[12pt,a4paper]{article}

\usepackage{setspace}
\usepackage{amsmath}
\usepackage{amsthm}
\usepackage{amsfonts}
\usepackage{mathtools}
\usepackage{graphicx}
\usepackage{caption}
\usepackage{subcaption}
\usepackage{xcolor}

\newtheorem{theorem}{Theorem}
\newtheorem{lemma}{Lemma}

\numberwithin{theorem}{section}
\numberwithin{Definition}{section}
\numberwithin{lemma}{section}
\numberwithin{Algorithm}{section}
\numberwithin{equation}{section}
\usepackage{hyperref}
\hypersetup{
    colorlinks=true,
    linkcolor=blue,
    filecolor=magenta,      
    urlcolor=cyan,
}

\usepackage{amsmath}
\usepackage{amsfonts}

\usepackage{longtable}
\usepackage{algorithm}

\usepackage[font=small,labelfont=bf]{caption}

\usepackage[]{algpseudocode}
    \usepackage{multirow}
\usepackage{geometry}
\usepackage{setspace}
\usepackage{float}
\usepackage[misc]{ifsym}

\makeatletter
\def\@cline#1-#2\@nil{%
  \omit
  \@multicnt#1%
  \advance\@multispan\m@ne
  \ifnum\@multicnt=\@ne\@firstofone{&\omit}\fi
  \@multicnt#2%
  \advance\@multicnt-#1%
  \advance\@multispan\@ne
  \leaders\hrule\@height\arrayrulewidth\hfill
  \cr
  \noalign{\nobreak\vskip-\arrayrulewidth}}
\makeatother

\begin{document}
	\title{Dynamic Non-Diagonal Regularization in Interior Point Methods for Linear and Convex Quadratic Programming}
\author{Spyridon Pougkakiotis$^1$   \and  Jacek Gondzio$^1$}

\maketitle

\noindent \textbf{Abstract}\\
\noindent In this paper, we present a dynamic non-diagonal regularization for interior point methods. The non-diagonal aspect of this regularization is implicit, since all the off-diagonal elements of the regularization matrices are cancelled out by those elements present in the Newton system, which do not contribute important information in the computation of the Newton direction. Such a regularization has multiple goals. The obvious one is to improve the spectral properties of the Newton system solved at each iteration of the interior point method. On the other hand, the regularization matrices introduce sparsity to the aforementioned linear system, allowing for more efficient factorizations. We also propose a rule for tuning the regularization dynamically based on the properties of the problem, such that sufficiently large eigenvalues of the non-regularized system are perturbed insignificantly. This alleviates the need of finding specific regularization values through experimentation, which is the most common approach in literature. We provide perturbation bounds for the eigenvalues of the non-regularized system matrix and then discuss the spectral properties of the regularized matrix. Finally, we demonstrate the efficiency of the method applied to solve standard small and medium-scale linear and convex quadratic programming test problems.


\section{Introduction} \label{section: intro}

\noindent In this paper, we are concerned with finding the solution of linear and convex quadratic programming problems, using an infeasible primal-dual interior point method. Such methods are called infeasible due to the fact that they allow intermediate iterates, produced by the algorithm, to be infeasible for the problem under consideration. They are called primal-dual, because they operate on both the primal and the dual space. Interior Point Methods (IPMs) deal with the inequality constraints of the problem by introducing logarithmic barriers in the objective, which penalize when any of the inequality constraints is close to being violated. At each iteration, the optimality conditions of the barrier problems are formed and one (or a few) steps of Newton method are applied to them. There is vast available literature on interior point methods and we refer the interested reader to \cite{paper_16} for an extended literature review. \\

\noindent Most implementations transform the Newton system into a symmetric indefinite system of linear equations, which when solved, determines the Newton direction. The latter constitutes the main computational effort and challenge for IPMs. At every iteration of the method, the system matrix as well as the right hand side change. There are three main reasons indicating why solving such a system can be challenging. The most obvious one, is that the dimension of such systems can be very large, which makes the task of solving them expensive in terms of processing time and memory requirements. A second important challenge, inherent in interior point methods, is that as the algorithm approaches optimality, the systems that we have to solve become increasingly ill-conditioned. Finally, a rank deficient constraint matrix can result in a singular Newton system matrix. It is well known that the latter two difficulties can be addressed by the use of some regularization technique, at the expense of solving a perturbed problem, \cite{paper_15}. \\

\noindent Such regularization techniques, embedded in the interior-point framework for solving linear and convex quadratic programming problems, have been previously proposed in the literature. For example, in \cite{paper_1}, a dynamic primal-dual regularization for interior point methods was derived. The authors solve a slightly altered symmetric indefinite system, to which a diagonal perturbation (regularization) has been introduced. This perturbation transforms the symmetric indefinite matrix into a quasi-definite one. It is proved in \cite{paper_4}, that such matrices are strongly factorizable. Hence, the regularized system can be factorized efficiently. The authors interpreted these regularization matrices as adding proximal terms to the primal and dual objective functions. The values of these perturbations are chosen dynamically during the factorization of the system matrix, where potentially unstable pivots are regularized stronger (using some pre-specified ``large" regularization value), while safer ones are almost not regularized at all. In \cite{paper_2}, based on this proximal point interpretation given in \cite{paper_1}, the authors proposed a primal-dual pair of regularized models, where the duality correspondence arises by setting the regularization variables as proximal terms. They observed that for specific parameter values, this primal-dual regularized model is \textit{exact}, that is it yields an optimal solution which is also an optimal solution of the respective non-regularized primal-dual pair. There, the authors introduced two uniform diagonal regularization matrices whose values were tuned experimentally over a variety of problems. A similar regularization was also used in \cite{paper_5}. It is worth mentioning that similar ideas have also been applied in IPMs suitable for general non-linear optimization problems (see \cite{paper_17,paper_18}). \\

\noindent In this paper, we are taking a different approach. We observe that when an IPM progresses and approaches optimality, significant part of the primal-dual variables approaches zero fast and hence becomes negligible. Yet it is not straightforward how the algorithm might exploit this feature. The proposed method attempts to do so. The method dynamically chooses a suitable regularization for the symmetric indefinite system and effectively ``annihilates" the effects of those parts of it, which do not contribute important information to the computation of the Newton direction. The proposed technique involves non-diagonal regularization matrices. However, their non-diagonal terms are only implicit; they do not need to be computed because they are immediately cancelled by other terms present in the linear system. Hence, the effect of adding such non-diagonal regularization is making the Newton system more sparse and therefore easier. In contrast to other previously developed approaches, this regularization is dynamically tuned based on the problem properties. We develop an approach which attempts to capture the needs of an arbitrary problem and regularize its system matrix accordingly. This alleviates the problem of finding specific regularization values that work well over a variety of problems. In general, the proposed approach is very conservative and regularizes the system as little as possible, while ensuring numerical stability. \\

\noindent The rest of the paper is organized as follows. In Section \ref{section: exact reg}, we summarize our notation and present the adopted model, based on which, we define our regularization matrices, firstly for linear and then for convex quadratic programming problems. For both cases, we provide arguments indicating why the proposed dynamic tuning of the regularization matrices is expected to introduce a controlled perturbation to the problem. In Section \ref{section: Spectral analysis}, we provide a spectral analysis, which shows the effect of the proposed regularization and gives specific bounds for the eigenvalues of the regularized system matrix. In Section \ref{section: Implementation}, we provide the algorithmic scheme along with some implementation details and numerical results, and finally in Section \ref{section: conlcusions} we derive our conclusions.

\section{Exact Primal-Dual Regularization} \label{section: exact reg}
\subsection{Notation}
\noindent Given an arbitrary symmetric square matrix $Q$, we denote positive semi-definiteness (positive definiteness) by $Q \succeq 0$ ($Q \succ 0$). We denote the Euclidean norm (2-norm) as $\|\cdot \|$. Any other norm, will be specified by a subscript. For example, the $\infty$-norm, is denoted as $\| \cdot \|_{\infty}$. We denote by $e$ the column vector of ones of appropriate dimension. Given a set of indices, say $\mathcal{B}$, $e_{\mathcal{B}}$ denotes the vector of ones with dimension equal to the cardinality of $\mathcal{B}$, that is: $e_{\mathcal{B}} \in \mathbb{R}^{|\mathcal{B}|}$. For an arbitrary matrix, say $A$, $A_{\mathcal{B}}$ denotes the sub-matrix whose columns and rows are indicated from the set of indices $\mathcal{B}$. Similarly, $A_{\mathcal{B}\mathcal{N}}$ contains rows of $A$ that belong in $\mathcal{B}$ and columns of $A$ that belong in $\mathcal{N}$. Iterates of the algorithm are denoted as $w_k = (x_k,r_k,s_k,y_k,z_k)$, where $k \in \mathbb{N}$ is the iteration counter. An optimal solution of the problem, is denoted as $w^* = (x^*,r^*,s^*,y^*,z^*)$. Given a vector $x \in \mathbb{R}^n$, we denote by $X \in \mathbb{R}^{n \times n}$ the diagonal matrix that contains $x$ in its diagonal. To simplify the notation, if a matrix, say $X_k$, depends on the iteration $k$, we will omit the sub-script and state the dependence whenever it is not obvious. When an arbitrary function, say $f$, depends on some parameter, say $\eta$, we denote this relation as: $f_{\eta}(\cdot)$. Given an arbitrary square matrix $B$, $\text{off}(B)$ denotes the square matrix that has the same off-diagonal elements as $B$ and has zeros in its diagonal. Similarly, $\text{diag}(B) = B - \text{off}(B)$. The j-th diagonal element of a square matrix $B$ will be denoted as: $(B)_{jj}$. $B^H$ denotes the conjugate (Hermitian) transpose of matrix $B$. We denote the smallest (largest) eigenvalue of an arbitrary matrix $B$, by $\lambda_{\min}(B)$ ($\lambda_{\max}(B)$). Similarly, the smallest (largest) singular value of an arbitrary matrix $B$ is denoted by $\sigma_{\min}(B)$ ($\sigma_{\max}(B)$). Finally, the set of all eigenvalues (spectrum) of an arbitrary matrix $B$, is denoted as $\lambda(B)$.
\subsection{Problem Formulation} \label{subsection: model}

\noindent We consider the following primal-dual pair of convex quadratic programming problems in the standard form:
\begin{equation} \label{non-regularized primal}
\tag{P}
\text{min}_{x} \ \big( c^Tx + \frac{1}{2}x^T Q x  \big), \ \ \text{s.t.}  \  Ax = b,    \ x \geq 0, 
\end{equation}
\begin{equation} \label{non-regularized dual}
\tag{D}
\text{max}_{x,y,z}  \ \big(b^Ty - \frac{1}{2}x^T Q x \big), \ \ \text{s.t.}\  -Qx + A^Ty + z = c,\ z\geq 0,
\end{equation}
\noindent where $c,\ x,\ z \in \mathbb{R}^n$, $b,\ y \in \mathbb{R}^m$, $A \in \mathbb{R}^{m\times n}$, $Q \succeq 0 \in \mathbb{R}^{n \times n}$. Without loss of generality, we assume that $m \leq n$. Note that if $Q = 0$, (\ref{non-regularized primal})-(\ref{non-regularized dual}) is a primal-dual pair of linear programming problems. If the problems under consideration are feasible, it can easily be verified that there exists an optimal primal-dual triple $(x,y,z)$ satisfying the Karush--Kuhn--Tucker (KKT) optimality conditions for this primal-dual pair (see for example Prop. 2.3.4 in \cite{book_2}).\\

\noindent Our model is based on the developments in \cite{paper_2}, \cite{paper_1} and \cite{paper_5}. More specifically, by applying a generalized primal-dual proximal point method on (\ref{non-regularized primal}), as in (\cite{paper_2,paper_18}), one can get the following pair of primal-dual regularized problems:
\begin{equation} \label{Primal Problem}
\tag{$P_r$}
\begin{split}
\text{min}_{x,r} &\ \bigg( c^Tx + \frac{1}{2}x^T Q x + \frac{1}{2} (r + \tilde{y})^T R_d(r+\tilde{y}) + \frac{1}{2}(x-\tilde{x})^T R_p (x-\tilde{x}) \bigg) \\  
\text{s.t.}  &\  Ax + R_d r = b, \ x \geq 0 ,  
\end{split}
\end{equation}

\begin{equation} \label{Dual Problem}
\tag{$D_r$}
\begin{split}
\text{max}_{x,y,z,s}&  \ \bigg( b^Ty - \frac{1}{2}x^T Q x - \frac{1}{2} (y-\tilde{y})^T R_d (y-\tilde{y}) - \frac{1}{2}(s+\tilde{x})^T R_p (s+\tilde{x}) \bigg) \\ 
\text{s.t.} &\ -Qx -R_p s+ A^Ty + z = c,\ z\geq 0,
\end{split}
\end{equation}
\noindent where $s \in \mathbb{R}^n$, $r \in \mathbb{R}^m$ are auxiliary variables introduced from the primal-dual application of the proximal point method and, $\ R_p \succeq 0 \in \mathbb{R}^{n \times n}$, $R_d \succeq 0\in \mathbb{R}^{m \times m}$ are the primal and dual regularization matrices respectively, that will be specified later. The duality correspondence follows after taking $r = y-\tilde{y}$ and $s = x - \tilde{x}$, where $\tilde{y}$ and $\tilde{x}$ are estimates of the dual and primal solutions $y^*,\ x^*$ respectively. Of course $R_p = 0,\ R_d = 0$ recovers the initial pair (\ref{non-regularized primal})-(\ref{non-regularized dual}). In \cite{paper_2}, the authors observe that this pair of regularized problems is \textit{exact} under some conditions on the estimates $\tilde{x},\ \tilde{y}$. In such a case, an optimal solution of (\ref{Primal Problem})-(\ref{Dual Problem}) is also an optimal solution of (\ref{non-regularized primal})-(\ref{non-regularized dual}). For more information about exactness of regularization, we refer the interested reader to \cite{paper_19}.\\

\noindent In \cite{paper_2,paper_17,paper_18}, models similar to (\ref{Primal Problem})-(\ref{Dual Problem}) are used, restricted however in the case where $R_p = \rho I$ and $R_d = \delta I$, for some positive values $\delta,\ \rho$. It is a well known fact, proved for the first time in \cite{paper_13}, that these regularization schemes can be interpreted as the primal-dual application of the standard proximal point method. However, our model does not specify the structure of the regularization matrices $R_p,\ R_d$. The only requirement is that these matrices are positive definite. As we commented previously, this model can be interpreted as the application of a generalized primal-dual proximal point method. Such methods, instead of adding the typical 2-norm in the objective function, make use of the so called D-functions. In fact, one could easily verify that any elliptic norm (defined by an arbitrary positive definite matrix) satisfies the conditions, given in \cite{paper_20,paper_14}, for being a D-function. In other words, our algorithm adds an elliptic norm in the objective, instead of the typical 2-norm. The focus of the paper however, prevents us from going deeper into these matters. For more about proximal point methods, we refer the reader to \cite{paper_21,paper_12,paper_13,paper_14,paper_20}, and the references therein.
\subsection{The Newton System} 

\noindent In order to solve the problems presented in the previous sub-section, using interior point methods, we proceed by replacing the non-negativity constraints with logarithmic barriers in the objective. In view of the previous, we obtain the following primal-dual regularized barrier problems:
\begin{equation} \label{Primal Barrier Problem}
\begin{split}
\text{min}_{x,r} &\ \bigg( c^Tx + \frac{1}{2}x^T Q x + \frac{1}{2} (r + \tilde{y})^T R_d(r+\tilde{y}) + \frac{1}{2}(x-\tilde{x})^T R_p (x-\tilde{x}) - \mu \sum_{j=1}^n \ln (x_j) \bigg)\\  
\text{s.t.}  &\  Ax + R_d r = b,  \\ 
\end{split}
\end{equation}
\begin{equation} \label{Dual Barrier Problem}
\begin{split}
\text{max}_{x,y,z,s}&  \bigg(\ b^Ty - \frac{1}{2}x^T Q x - \frac{1}{2} (y-\tilde{y})^T R_d (y-\tilde{y}) - \frac{1}{2}(s+\tilde{x})^T R_p (s+\tilde{x}) - \mu \sum_{j=1}^n \ln (z_j) \bigg)\\ 
\text{s.t.}\ & -Qx -R_p s+ A^Ty + z = c,\\
\end{split}
\end{equation}
\noindent in which non-negativity constraints $x>0$ and $z>0$ are implicit.\\
 
\noindent Forming the Lagrangian of the primal barrier problem, we get:
\begin{equation} \label{Lagrangian}
\begin{split}
\mathcal{L}_{\tilde{x},\tilde{y},\mu}(x,y,r) = &\  c^T x + \frac{1}{2}x^T Q x +  \frac{1}{2} (r + \tilde{y})^T R_d(r+\tilde{y})+  \\& + \frac{1}{2}(x-\tilde{x})^T R_p (x-\tilde{x}) - y^T(Ax + R_d r - b) - \mu \sum_{j=1}^n \ln (x_j).
\end{split}
\end{equation}
\noindent Now, we can form the first order optimality conditions of the problems by taking the gradient of (\ref{Lagrangian}) and equating it to zero, giving us the following block equations: 
\begin{equation*}
\begin{split}
\nabla_x \mathcal{L}_{\tilde{x},\tilde{y},\mu}(x,y,r) =&\ c + Qx + R_p (x-\tilde{x}) - A^T y - \mu X^{-1} e = 0,\\
\nabla_y \mathcal{L}_{\tilde{x},\tilde{y},\mu}(x,y,r) =&\ Ax + R_d r - b = 0,\\
\nabla_r \mathcal{L}_{\tilde{x},\tilde{y},\mu}(x,y,r) = &\ R_d(r + \tilde{y}) - R_d y = 0.\\
\end{split}
\end{equation*}
\noindent By looking at the optimality conditions of the dual barrier problem, we see that the final two conditions are:
\begin{equation*}
\begin{split}
R_p x - R_p (s + \tilde{x}) &= 0,\\
XZe &= \mu e.
\end{split}
\end{equation*}
\noindent We write the optimality conditions in the form of a function, $F_{\tilde{x},\tilde{y},\mu}(w)\ :\ \mathbb{R}^{3n+ 2m} \rightarrow \mathbb{R}^{3n+2m}$ and we want to solve:
\begin{equation} \label{optim. contitions}
F_{\tilde{x},\tilde{y}, \mu}(w) = 
\begin{bmatrix} 
c + Qx + R_p s - A^T y - z\\
R_d(r + \tilde{y}) - R_d y \\
R_p x - R_p (s + \tilde{x})\\
Ax + R_d r - b\\
XZe
\end{bmatrix}
= 
\begin{bmatrix}
0\\
0\\
0\\
0\\
\sigma \mu e
\end{bmatrix},
\end{equation}
\noindent  at each IPM iteration, where $w = (x,r,s,y,z),\ \mu >0$ is the barrier parameter and $\sigma \in\ ]0,1[$ is a centring parameter (which determines how fast $\mu$ is forced to decrease). We want to force $\mu \rightarrow 0$, since then, the solution of this system leads to the solution of (\ref{Primal Problem})-(\ref{Dual Problem}). Notice that (\ref{Primal Problem})-(\ref{Dual Problem}) is parametrized by the estimates $\tilde{x}$ and $\tilde{y}$. As observed in \cite{paper_2}, if these estimates are close enough to some optimal solution of (\ref{non-regularized primal})-(\ref{non-regularized dual}), then an optimal solution of (\ref{Primal Problem})-(\ref{Dual Problem}) is also an optimal solution of (\ref{non-regularized primal})-(\ref{non-regularized dual}). At the beginning of the $k$-th iteration of the IPM, we have available the iterate $w_k = (x_k,r_k,s_k,y_k,z_k)$, the barrier parameter $\mu_k = \frac{x_k^Tz_k}{n}$ and we choose a value for the centring parameter $\sigma_{k} \in\ ]0,1[$. Following the developments in \cite{paper_1,paper_2,paper_17,paper_18}, for proximal point methods, we update the estimates of $x^*,\ y^*$ as $\tilde{x} = x_k,\ \tilde{y} = y_k$. Next, Newton method is applied to the mildly non-linear system (\ref{optim. contitions}). After evaluating the Jacobian of $F_{\tilde{x},\tilde{y},\mu}(w)$, the Newton direction is determined at each IPM iteration by solving a system of the following form:
\begin{equation} \label{Newton System}
\begin{bmatrix} 
Q & 0 & R_p & - A^T & -I\\
0 & R_d & 0 & -R_d & 0\\
R_p & 0 & -R_p & 0 & 0\\
A & R_d & 0 & 0 & 0\\
Z &  0 & 0 & 0  &X 
\end{bmatrix}
\begin{bmatrix}
\Delta x\\
\Delta r\\ 
\Delta s\\
\Delta y\\
\Delta z
\end{bmatrix}
= -
\begin{bmatrix}
c + Qx_k + R_p s_k - A^T y_k - z_k\\
R_d(r_k + \tilde{y}) - R_d y_k\\
R_p x_k - R_p (s_k + \tilde{x})\\
Ax_k + R_d r_k - b\\
 XZe - \sigma_{k} \mu_k e
\end{bmatrix}.
\end{equation}
\noindent Notice that the matrices $X, Z, R_p$ and $R_d$ all depend on the iteration $k$ of the algorithm. Once the Newton direction $\Delta w = (\Delta x,\Delta r,\Delta s,\Delta y,\Delta z)$ is computed, the algorithm chooses a step-length $a_k \in\ ]0,1]$ and sets the new iterate to $w_{k+1} = w_k + a_k \Delta w$. In order to compute the Newton direction efficiently, we want to eliminate some variables of $(\ref{Newton System})$. Since we set $\tilde{y} = y_k$, the second block equation of (\ref{Newton System}) gives:
\[R_d \Delta r - R_d \Delta y = -R_d(r_k+y_k) + R_d y_{k}, \]
\noindent and if $R_d \succ 0$, we get the following relation:
\begin{equation} \label{Delta r recovery}
\Delta y = r_k + \Delta r,
\end{equation}
\noindent Similarly, by looking at the third block equation of (\ref{Newton System}) and substituting $\tilde{x} = x_{k}$, we get:
\[R_p \Delta x - R_d \Delta s = R_p (s_k + x_k) - R_p x_{k},\]
\noindent and if $R_p \succ 0$ we have that:
\begin{equation} \label{Delta s recovery}
\Delta x = s_k + \Delta s.
\end{equation}
\noindent Note that we always use either $R_d \succ 0$ or $R_d = 0$ and similarly, either $R_p \succ 0$ or $R_p = 0$. Hence, the previous two relations are either well-defined or absent. Using (\ref{Delta r recovery}) and (\ref{Delta s recovery}) to eliminate $\Delta r$ and $\Delta s$, we can reduce (\ref{Newton System}) to the following system:

\begin{equation} \label{3x3 Newton system}
\begin{bmatrix} 
-(Q + R_p) & \quad  A^T &\quad I\\
A & R_d & 0\\
Z &  0   &X 
\end{bmatrix}
\begin{bmatrix}
\Delta x\\ 
\Delta y\\
\Delta z
\end{bmatrix}
= 
\begin{bmatrix}
c + Qx_k - A^T y_k - z_k\\
b-Ax_k\\
\sigma_{k} \mu_k e - XZe
\end{bmatrix}.
\end{equation}

\noindent Next, we proceed by eliminating $\Delta z$. For that purpose, we have from the third row of (\ref{3x3 Newton system}) that:
\begin{equation} \label{Delta z recovery}
\Delta z = -X^{-1}Z \Delta x - Ze + \sigma_{k} \mu_k X^{-1}  e. 
\end{equation}
\noindent Substituting (\ref{Delta z recovery}) into the first row of (\ref{3x3 Newton system}), we get the following reduced symmetric system (so called \textit{Augmented System}):
\begin{equation} \label{Augmented System}
\begin{bmatrix} 
-(Q+\Theta^{-1} + R_p) &\quad   A^T \\
A & R_d
\end{bmatrix}
\begin{bmatrix}
\Delta x\\ 
\Delta y
\end{bmatrix}
= 
\begin{bmatrix}
c + Qx_k - A^T y_k - \sigma_{k} \mu_k X^{-1} e\\
b-Ax_k
\end{bmatrix},
\end{equation}
\noindent where $\Theta = X Z^{-1}$. In the case of linear programming ($Q = 0$) or when solving quadratic separable problems (in which case $Q$ is diagonal), it may be beneficial to further eliminate $\Delta x$ from (\ref{Augmented System}), which will end up at the so called \textit{normal equations}. However, one should note that this is not a good idea when it comes to general convex quadratic programming problems, since pivoting on the (1,1) block of (\ref{Augmented System}) could result in a dense system, even in cases where both $A$ and $Q$ are sparse. Having said that, we can eliminate $\Delta x$ by looking at the first block equation of (\ref{Augmented System}), which gives:
\begin{equation}\label{Delta x recovery}
\Delta x = (Q +\Theta^{-1}+R_p)^{-1} A^T \Delta y - (Q+\Theta^{-1}+R_p)^{-1} (c+Qx_k-A^T y_k - \sigma_{k} \mu_k X^{-1}e),
\end{equation}
\noindent and by substituting (\ref{Delta x recovery}) into the second row of (\ref{Augmented System}), we get the normal equations:
\begin{equation} \label{Normal Equations}
\big[A(Q+\Theta^{-1}+R_p)^{-1} A^T + R_d \big]\Delta y = \xi, 
\end{equation}
\noindent where 
\[ \xi = b - Ax_k + A(Q +\Theta^{-1} + R_p)^{-1} ( c+Qx_k-A^T y_k - \sigma_{k} \mu_k X^{-1}e),\]

\noindent in which the system matrix is symmetric and positive definite.\\

\noindent The proposed model differs from the one derived in \cite{paper_2} in that it allows the use of general positive definite regularization matrices. For example, if $R_p,\ R_d$ are non-diagonal matrices, then this would amount to the primal and dual application of a generalized proximal point method that adds an elliptic norm in the objective, instead of the typical 2-norm that is employed in standard proximal point methods. Notice that at every iteration of the algorithm, $R_p,\ R_d,\ \tilde{x}$ and $\tilde{y}$ are updated. In other words, (\ref{Primal Problem})-(\ref{Dual Problem}) represents a sequence of sub-problems. At every such sub-problem, we apply a single iteration of the interior point method. How $R_p$ and $R_d$ are updated will be presented in the following sub-section.

\subsection{The Regularization Matrices} \label{subsection: the regularization matrices}

\noindent As IPM approaches optimality, the diagonal matrix $\Theta$ contains elements that converge to zero and others that diverge to infinity. This is because $\mu_k \rightarrow 0$ and we force the complementarity conditions to be approximately satisfied ($XZe \approx \sigma_{k} \mu_k e$) . As a consequence, the matrices in (\ref{Augmented System}) and (\ref{Normal Equations}) become extremely ill-conditioned. On top of that, it is often the case due to modelling choices, that the constraint matrix $A$ is not of full row rank, which makes the system matrices singular. It is well known, as shown by Armand and Benoist \cite{paper_15}, that both these problems can be addressed with the use of regularization. The most common approach in the literature, is the addition of two diagonal regularization matrices, say $R_p,\ R_d$, whose values are tuned experimentally over a variety of problems (\cite{paper_5,paper_1,paper_2,paper_15}). \\

\noindent Roughly speaking, the goals of a regularization method for IPMs are (\cite{paper_15,paper_1,paper_5,paper_22,paper_17,paper_18}):
\begin{enumerate}
\item to improve the spectral properties of the matrices in (\ref{Augmented System}) and (\ref{Normal Equations}),
\item without significantly perturbing the previous systems,
\item while preserving the sparsity of the problem and the computational efficiency of the method.\\
\end{enumerate}
\noindent To the best of our knowledge, most of the regularization methods in literature manage to achieve the first and the third regularization goals, failing however to achieve the second goal with certainty. This is the case since these regularization methods are tuned experimentally. Hence, they do not rely on the properties of the problem itself, and as a consequence, such regularization values can only be good for some problems and poor for others. The proposed method takes a different approach, by introducing two non-diagonal regularization matrices $R_p$ and $R_d$, which are tuned based on the properties of the problem. Of course one could argue that this may disturb the sparsity and as a consequence the computational efficiency of the method, however, these non-diagonal matrices are created implicitly. As we will show later, not only the sparsity is preserved, but in fact it is improved. \\

\noindent As we already mentioned, as IPM approaches optimality, the matrix $\Theta$ contains some very large and some very small elements. The proposed regularization exploits this inherent feature of the method and splits the columns of the problem matrix in two sets, say $\mathcal{N}$ and $\mathcal{B}$ such that:
\begin{equation*}
\begin{split}
\forall\ j \in \mathcal{N}\ :\ x_j \rightarrow 0 ,\ z_j \rightarrow \hat{z}_j > 0 \Rightarrow\ & \ (\Theta)_{jj} = \frac{x_j}{z_j} \approx \frac{x_j z_j}{z_j^2} = O(\mu)\\
\forall\ j \in \mathcal{B}\ :\ x_j \rightarrow \hat{x}_j > 0,\ z_j \rightarrow 0  \Rightarrow\ & \ (\Theta)_{jj}= \frac{x_j}{z_j} \approx \frac{x_j^2}{x_j z_j} = O(\mu^{-1}),
\end{split}
\end{equation*}
\noindent where $|\mathcal{N}| = n_1$ and $| \mathcal{B}| = n_2$, with $n_1 + n_2 = n$.  Notice that the previous splitting captures all the columns only if the method converges to a strictly complementary solution (that is the limit point satisfies: $\hat{x}^T\hat{z} = 0$ and $\hat{x}_j + \hat{z}_j > 0,\ \forall\ j$). In the quadratic programming case, a strictly complementary solution may not exist. Hence, there might exist some indices $j \subseteq \{1,\cdots,n\}$ for which: $x_j \rightarrow 0$ and $z_j \rightarrow 0$. In such a case, it is unknown whether the value of $\Theta_{jj}$ will be small or large. We can assume, without loss of generality, that any such indices will be classified as elements of $\mathcal{B}$ (although in practice this would depend on the value of $\Theta_{jj}$, as we will show later). Of course for the case of linear programming ($Q = 0$), it is a well-known fact (see for example \cite{paper_26}) that a strictly complementary solution always exists, if the problems are feasible. Moreover, as shown in \cite{paper_23,paper_24}, primal-dual IPMs converge to such an optimal solution. If a strictly complementary solution exists for the quadratic programming case, it is shown in \cite{paper_25}, that an infeasible primal-dual IPM which reduces the constraints violation at the same rate as $\mu$ is reduced, produces iterates that converge to a strictly complementary solution.\\

\noindent In what follows, we present the construction of the regularization for the case of linear programming and then we suggest an extension for convex quadratic programming. 
\subsubsection{Linear Programming} \label{subsubsection: regularization matrices LP}

\noindent For the case of linear programming we employ a dual regularization, that is, in (\ref{Newton System}) we set $R_p = 0$ and only use $R_d \succ 0$ to improve the spectral properties of the problem. Given this set-up, and by permuting the columns so that the first $n_1$ of them correspond to indices in $\mathcal{N}$ while the remaining correspond to indices in $\mathcal{B}$, the augmented system in (\ref{Augmented System}) takes the form:
\begin{equation}  \label{LP permuted augmented system}
\begin{bmatrix} 
-\Theta_{\mathcal{N}}^{-1} &  0 & A_{\mathcal{N}}^T \\
0 &  -\Theta_{\mathcal{B}}^{-1}& A_{\mathcal{B}}^T\\
A_{\mathcal{N}} & A_{\mathcal{B}} & R_d
\end{bmatrix}
\begin{bmatrix}
\Delta x_{\mathcal{N}}\\
\Delta x_{\mathcal{B}} \\ 
\Delta y
\end{bmatrix}
= 
\begin{bmatrix}
c_{\mathcal{N}} - A_{\mathcal{N}}^T y_k - \sigma_{k} \mu_k X_{\mathcal{N}}^{-1} e_{\mathcal{N}}\\
c_{\mathcal{B}} - A_{\mathcal{B}}^T y_k - \sigma_{k} \mu_k X_{\mathcal{B}}^{-1} e_{\mathcal{B}}\\
b-Ax_k
\end{bmatrix},
\end{equation}
\noindent where $A_{\mathcal{N}} \in \mathbb{R}^{m \times n_1}$ and $A_{\mathcal{B}} \in \mathbb{R}^{m\times n_2}$. Pivoting on the first $n_1$ columns of (\ref{LP permuted augmented system}), gives the \textit{partially reduced augmented system}:
\begin{equation}  \label{LP partially reduced augmented system}
\begin{bmatrix} 
  -\Theta_{\mathcal{B}}^{-1}& A_{\mathcal{B}}^T\\
 A_{\mathcal{B}} & A_{\mathcal{N}} \Theta_{\mathcal{N}} A_{\mathcal{N}}^T + R_d
\end{bmatrix}
\begin{bmatrix}
\Delta x_{\mathcal{B}} \\ 
\Delta y
\end{bmatrix}
= 
\begin{bmatrix}
c_{\mathcal{B}} - A_{\mathcal{B}}^T y_k - \sigma_k \mu_k X_{\mathcal{B}}^{-1} e_{\mathcal{B}}\\
b-Ax_k + A_{\mathcal{N}} \Theta_{\mathcal{N}} ( c_{\mathcal{N}}-A_{\mathcal{N}}^T y_k - \sigma_{k} \mu_k X_{\mathcal{N}}^{-1}e_{\mathcal{N}})
\end{bmatrix}.
\end{equation}
\noindent Since we know that $\Theta_{\mathcal{N}} \rightarrow 0$, we expect that the magnitude of  $\|A_{\mathcal{N}} \Theta_{\mathcal{N}} A_{\mathcal{N}}^T\|$ will be small when the method approaches optimality. Intuitively, our goal is to create a regularization matrix that will implicitly absorb the off-diagonal elements of $A_{\mathcal{N}} \Theta_{\mathcal{N}} A_{\mathcal{N}}^T$ (promoting sparsity) and regularize the system with values having a slightly larger magnitude to that of the elements which were absorbed. For this class of problems, we will focus on solving the normal equations. Given (\ref{LP partially reduced augmented system}), we can form the normal equations by eliminating $\Delta x_{\mathcal{B}}$, which gives the following system:
$$ \big[A_{\mathcal{B}}\Theta_{\mathcal{B}} A_{\mathcal{B}}^T + A_{\mathcal{N}} \Theta_{\mathcal{N}} A_{\mathcal{N}}^T + R_d\big]\Delta y= b - Ax_k + A\Theta ( c-A^T y_k - \sigma_{k} \mu_k X^{-1}e).$$
\noindent  We choose the following dual regularization matrix:
\begin{equation}\label{dual regularization matrix}
R_d =  \big(\Delta_d - \text{off}(A_{\mathcal{N}}\Theta_{\mathcal{N}} A_{\mathcal{N}}^T)\big), 
\end{equation} 
\noindent where $\Delta_d$ is a diagonal matrix chosen such that $R_d \succ 0$ and diagonally dominant, that is:
\[(\Delta_{d})_{ii} > \sum_{j=1,j \neq i}^{m} |(A_{\mathcal{N}}\Theta_{\mathcal{N}} A_{\mathcal{N}}^T)_{ij} |,\ \ \forall \ i = 1,\cdots,m.\]
\noindent For computational efficiency and numerical stability, we choose $\Delta_d = \delta_{d,k} I_m$, with:
\begin{equation} \label{LP Delta_d definition}
\delta_{d,k} = (\max_{j}(\Theta_{\mathcal{N}})_{jj}) \|A_{\mathcal{N}} A_{\mathcal{N}}^T\|_{\infty}.
\end{equation}
\noindent Observe that the regularization matrix given in (\ref{dual regularization matrix}), strongly depends on the properties of the problem as well as on the iteration $k$ of the IPM. In order to control which elements enter the set $\mathcal{N}$, at every iteration $k$, we enforce the following condition:
\begin{equation} \label{LP Condition for entering N}
\max_j (\Theta_{\mathcal{N}})_{jj} \| A A^T\|_{\infty} \leq \text{reg}_{thr,k},
\end{equation}
\noindent where $\text{reg}_{thr,k}$ is set to $1$ at the beginning of the optimization ($k = 0$), and is decreased at the same rate as $\mu_k$ (i.e. $\text{reg}_{thr,k} = O(\mu_k)$). Once $\text{reg}_{thr,K}$ becomes smaller than a predefined value, say $\epsilon > 0$, for some large $K \geq 1$, we fix it to this value ($\text{reg}_{thr,k} = \epsilon,\ \forall\ k \geq K$). The choice of $\epsilon$ will be specified later. Note that (\ref{LP Condition for entering N}) ensures that $\delta_{d,k} < \text{reg}_{thr,k}$, at every iteration. In order to show that sparsity is improved, we form again the normal equations' matrix using the definition of $R_d$ to get:
\[A\Theta A^T + R_d = A_{\mathcal{B}}\Theta_{\mathcal{B}} A_{\mathcal{B}}^T + \text{diag}\big(A_{\mathcal{N}} \Theta_{\mathcal{N}} A_{\mathcal{N}}^T \big)+ \Delta_d. \]
\noindent From the previous one can easily observe that the sparsity of the normal equations is improved, since some off-diagonal elements of the matrix have been absorbed by the regularization.\\ 

\noindent  Since $\text{reg}_{thr,k}$ is not allowed to go to zero as $\mu_k \rightarrow 0$, we would like to know how much we perturb the Newton system, by having it fixed to some value $\epsilon > 0$, when the method is close to optimality. In the rest of this subsection, we compute some perturbation bounds, which depend on the value of $\text{reg}_{thr}$.

\paragraph*{Motivation}
$\ $\\
\noindent Now that we have defined the regularization matrix for the case of linear programming problems, let us provide a motivation for this choice. Firstly, note that the proposed regularization has multiple objectives. On the one hand, we want to find a good criterion for tuning a uniform dual regularization matrix $\delta_{d,k} I$ based on the properties of the problem, such that the non-regularized problem matrix is not perturbed significantly while its spectral properties are improved. On the other hand, we use this uniform dual regularization value as a cut-off point, for dropping the smallest off-diagonal elements in the normal equations matrix, improving the computational efficiency of the method. In what follows we will provide an analysis indicating why the uniform dual regularization that we introduce is expected not to perturb the problem significantly. Then, we will show that further dropping the off-diagonal elements introduces a controlled perturbation.\\

\noindent Based on the previous, let us assume for now that $R_d = \delta_{d,k} I$, where $\delta_{d,k}$ is defined as in (\ref{LP Delta_d definition}). For simplicity of notation, we omit the iteration subscript in $\delta_d$ and we let:
\begin{equation*}
M = \begin{bmatrix}
-\Theta^{-1} & A^T\\
A & 0
\end{bmatrix},\ \ 
E = \begin{bmatrix}
0 & 0\\
0 &\quad \delta_d I
\end{bmatrix}.
\end{equation*}
\noindent We want to analyse the difference in the eigenvalues of the matrices $M$ and $M+E$. For the rest of this sub-section, let $\lambda_i$ denote the i-th smallest eigenvalue of $M$, $\tilde{\lambda}_i$ the i-th smallest eigenvalue of $M+E$, and $\lambda_i(t)$ the i-th smallest eigenvalue of $M+tE$, with $t \in [0,1]$. The smallest eigenvalues of $M$ (in the absolute value sense) will be increased after the addition of $E$ and this is of course desirable, since this was the main motivation for introducing the regularization. The following analysis provides perturbation bounds only for eigenvalues of $M$ that satisfy $|\lambda_i |> 2\|E\|$. We will assume also that the eigenvalues that we analyse are simple (i.e. their algebraic multiplicity is 1). The analysis can be extended to multiple eigenvalues, however it gets unnecessarily complicated. Such an analysis is derived in the appendix of \cite{paper_10}. Let us now state a lemma derived in \cite{book_1}.
\begin{lemma} \label{lemma partial derivative of eigenvalue}
\noindent Let $M$, $E$ be square Hermitian matrices. Denote by $\lambda_i(t)$ the i-th smallest eigenvalue of $M+tE$ and consider the eigenvector function $x(t)$ such that: $(M+tE)x(t) = \lambda_i(t)x(t)$, with $\|x(t)\| = 1$, for some $t \in [0,1]$. If $\lambda_i(t)$ is simple, then:
\[\frac{\partial \lambda_i(t)}{\partial t} = x(t)^H E x(t).\]
\end{lemma}
\noindent As observed in \cite{paper_10}, if the eigenvector $x(t)$ has small components in the positions corresponding to the dominant elements of $E$, then $\frac{\partial \lambda_i(t)}{\partial t}$ is expected to be small. Let us now provide the following lemma, based on the developments in \cite{paper_9}.
\begin{lemma}
Let $\lambda_i \neq 0$ be an eigenvalue of $M$ and $Mx = \lambda_i x$, with $\|x\| = 1$. Partitioning $x = [x_1^H\ x_2^H]^H$, we have:
\[\|x_2\| \leq \frac{\|A\|}{\sqrt{\lambda_i^2+\|A\|^2}}.\]
\end{lemma}
\noindent {\textbf{Proof}}
\noindent The proof follows exactly the developments in \cite{paper_9}, but we provide it here for completeness. From the second block equation of $M x = \lambda_i x$, we have:
\[A x_1 = \lambda_i x_2 \Rightarrow x_2 = \frac{1}{\lambda_i} Ax_1,\]
\noindent where the latter is well defined since we have assumed that $\lambda_i \neq 0$. By taking norms on both sides in the previous equation, we get:
\[\|x_2\| \leq \frac{1}{|\lambda_i|} \|A\| \|x_1\|.\]
\noindent But $\|x\| = 1 \Rightarrow \|x_1\| = \sqrt{1 - \|x_2\|^2}$. Hence, we have:
\[\|x_2\| \leq \frac{\|A\| \sqrt{1-\|x_2\|^2}}{|\lambda_i|}.\]
\noindent By solving the previous inequality, we get: 
\[\|x_2\| \leq \frac{\|A\|}{\sqrt{\lambda_i^2+\|A\|^2}},\]	
\noindent which completes the proof.\qed

\noindent The following lemma will be a useful tool for the analysis. We omit its trivial proof.
\begin{lemma} \label{lemma monotone increasing function}
\noindent Let $f(x) = \frac{x}{\sqrt{a+x^2}}$, where $a > 0$. Then, $f(x)$ is a monotone increasing function for $x > 0$.
\end{lemma}
\noindent Let us now bound the second block of the eigenvector function $x_2(t)$ based on the developments in \cite{paper_9}.
\begin{lemma} \label{lemma LP motivation x2(t)}
\noindent Assume that $\lambda_i \neq 0$ is the i-th smallest eigenvalue of $M$. Consider the eigenvector function $x(t)$ such that: $(M+ tE)x(t) = \lambda_i(t)x(t)$, with $\|x(t)\| = 1$, $\forall\ t \in [0,1]$. Partitioning $x(t) = [x_1(t)^H\ x_2(t)^H]^H$ and assuming that $|\lambda_i| > 2\|E\|$, we have that:
\[\|x_2(t)\| \leq \frac{\|A\|}{\sqrt{(|\lambda_i| - 2\|E\|)^2 + \|A\|^2}}.\]
\end{lemma}
\noindent {\textbf{Proof}}
\noindent We omit the proof which follows from Lemma \ref{lemma monotone increasing function} combined with the previous developments. The interested reader can view \cite{paper_9}, Lemma 2.8, for a detailed derivation which can directly be applied in our context.
\qed
\noindent Let us now derive the following theorem which bounds the difference between the i-th smallest eigenvalues of the matrices $M$ and $M+E$ respectively. 
\begin{theorem} \label{Theorem LP reg motivation}
\noindent Let $\lambda_i$ and $\tilde{\lambda}_i$ be the respective i-th smallest eigenvalues of $M$ and $M+E$ and define $\phi_i = \frac{\|A\|}{\sqrt{(|\lambda_i| - 2\|E\|)^2 + \|A\|^2}}$. For every $i$ such that $|\lambda_i| > 2\|E\|$ we have that:
\[ |\lambda_i - \tilde{\lambda}_i | \leq \|E\| \phi_i^2.\]
\end{theorem}
\noindent {\textbf{Proof}}
\noindent From Lemma \ref{lemma partial derivative of eigenvalue} and Lemma \ref{lemma LP motivation x2(t)} it follows that:
\begin{equation*}
\begin{split}
|\lambda_i - \tilde{\lambda}_i| =\ & |\lambda_i(0) - \lambda_i(1)| \\
 =\ & \bigg| \int_0^1 x(t)^H E x(t) dt \bigg| \\
 =\ & \bigg|\int_0^1 x_2(t)^H \delta_d I x_2(t) dt \bigg|\\
 =\ & \delta_d \int_0^1 \|x_2(t)\|^2 dt\\
 \leq\ & \|E\| \phi_i^2 = \delta_d \phi_i^2.
\end{split}
\end{equation*}
The proof is complete. \qed $\ $\\ 
\noindent Note that, since $\phi_i < 1$, the latter is a tighter bound than the general bound provided by Weyl's inequality, given that the eigenvalue under consideration is larger than $2\|E\|$. From the previous results we can draw several useful observations. As we already stated, the smaller the components of $x_2(t)$ are, the smaller $\frac{\partial \lambda_i(t)}{\partial t}$ is expected to be. But $x_2(t)$ is bounded by $\phi_i$. Hence, the smaller $\phi_i$ is, the more insensitive the eigenvalue $\lambda_i$ is to the perturbation $\|E\| = \delta_d$. In fact, in the previous theorem we proved that the error in the eigenvalue is bounded by $\|E\| \phi_i^2$.\\

\noindent Let us now examine the nature of $\phi_i$. Firstly, one can see that it depends on the norm of the constraint matrix $A$, and from Lemma \ref{lemma monotone increasing function} we can observe that it is monotone increasing with respect to the norm of $A$. What this tells us, is that the smaller the norm of the constraint matrix $A$ is, the more insensitive the eigenvalues of matrix $M$ are to the perturbation $E$. Of course the latter holds only for eigenvalues that are sufficiently larger than $2\|E\|$. On the other hand, from the definition of $\phi_i$, we can see that it is beneficial to have a small $\|E\|$, since then, most of the eigenvalues of $M$ are expected to satisfy: $|\lambda_i| > 2\|E\|$. \\

\noindent We now shift our attention to the proposed tuning of the regularization parameters. From (\ref{LP Condition for entering N}), the set of indices $\mathcal{N}$ is such that: $\max_j (\Theta_{\mathcal{N}})_{jj} \|AA^T\|_{\infty} \leq \text{reg}_{thr}$. Also, from (\ref{LP Delta_d definition}), we have that $\delta_d = \max_j (\Theta_{\mathcal{N}})_{jj} \|A_{\mathcal{N}} A_{\mathcal{N}}^T\|_{\infty}$. By combining the previous, we get:
\[\|E\| = \delta_d  \leq \frac{\text{reg}_{thr} \|A_{\mathcal{N}}A_{\mathcal{N}}^T\|_{\infty}}{\| A A^T \|_{\infty}}.\]

\noindent  Observe that, if $\|AA^T\|_{\infty}$ is large, we allow few columns to enter the partition $\mathcal{N}$.  In this case, $\phi_i$ is expected to be close to 1 for most of the eigenvalues  $\lambda(M)$. On the other hand, $|\mathcal{N}|$ is increased if the infinity norm of $AA^T$ is small, and in such a case, $\phi_i$ is expected to be small for many eigenvalues of the system matrix $M$. A more sophisticated choice for the regularization value based on the derived bounds is possible, however, the proposed regularization has two goals, that is not to perturb the system significantly while introducing sparsity to the problem, and hence the definition of $\delta_d$ is computationally advantageous for that. Note that taking advantage of the previously presented bounds indicates that the sufficiently large (in the absolute value sense) eigenvalues of the system matrix ($\gg 2\delta_d$) will be perturbed almost insignificantly. If some eigenvalues of the matrix are very small, the previous arguments break down. We will derive lower bounds for these eigenvalues in the next section.\\

\noindent Having introduced the diagonal uniform regularization $\delta_dI$, let us examine the effect of further dropping the off-diagonal elements $\text{off}(A_{\mathcal{N}} \Theta_{\mathcal{N}} A_{\mathcal{N}}^T)$ from the normal equations (\ref{Normal Equations}). For that, we define $K = A\Theta A^T + \delta_d I$ and $R = \text{off}(A_{\mathcal{N}} \Theta_{\mathcal{N}} A_{\mathcal{N}}^T)$ and consider the following generalized eigenvalue problem:
\begin{equation} \label{LP generalized eigenvalue problem}
u^TRu = \lambda u^T K u.
\end{equation}
\noindent The previous is well defined since $K \succ 0$. We will analyse the eigenvalues of $K^{-\frac{1}{2}}RK^{-\frac{1}{2}}$, which is similar to $K^{-1}R$. Now assume by contradiction that $\lambda_{\max}(K^{-\frac{1}{2}}RK^{-\frac{1}{2}}) \geq 1$. Then from (\ref{LP generalized eigenvalue problem}) and for some eigenvector $u$ corresponding to the maximum eigenvalue, we would have:
\[ u^T R u \geq u^T K u.\]
\noindent By adding $u^T \text{diag}(A_{\mathcal{N}} \Theta_{\mathcal{N}} A_{\mathcal{N}}^T)u$ to both sides of the previous inequality, we get:
\[0 \geq u^T (A_{\mathcal{B}} \Theta_{\mathcal{B}} A_{\mathcal{B}}^T)u +u^T \text{diag}(A_{\mathcal{N}} \Theta_{\mathcal{N}} A_{\mathcal{N}}^T)u + u^T \delta_d u,\]
\noindent which is a contradiction. Hence $\lambda_{\max}(K^{-\frac{1}{2}}RK^{-\frac{1}{2}}) < 1$. On the other hand, if we assume by contradiction that $\lambda_{\min}(K^{-\frac{1}{2}}RK^{-\frac{1}{2}} )\leq -1$, from (\ref{LP generalized eigenvalue problem}) and for an eigenvector $u$ corresponding to the minimum eigenvalue, we would get:
\begin{equation*}
\begin{split}
u^T R u \leq\ & -u^T Ku = -u^T(A\Theta A^T + \delta_dI)u \leq -\delta_d u^T u.\\
\end{split}
\end{equation*}
\noindent However, using (\ref{LP Delta_d definition}), we get $\delta_d + R \succ 0$, hence $-\delta_d u^T u < u^TRu$, which contradicts the previous inequality. Hence, $\lambda_{\min}(K^{-\frac{1}{2}}RK^{-\frac{1}{2}} )> -1$. Now, one can easily observe that: 
 $$K^{-1}(K-R) = I - K^{-1}R,\ \text{and}\ \rho(K^{-1}R) <1,$$
\noindent where $\rho(\cdot)$ is the spectral radius, and hence the eigenvalues of $K^{-1}(K-R)$ are clustered around 1. This supports the claim that further dropping the off-diagonal elements of the part of the normal equations corresponding to indices in $\mathcal{N}$, after adding a uniform dual regularization, introduces a controlled perturbation.
\subsubsection{Quadratic Programming} \label{subsubsection: regularization matrices QP}

Unlike the case of linear programming, for the case of quadratic programming we employ a primal-dual regularization, that is, we use both $R_p \succ 0$ and $R_d \succ 0$, as shown in (\ref{Newton System}), to improve the spectral properties of the problem. For this case, we modify the condition for allowing a column to enter the set $\mathcal{N}$, and at each iteration $k$, in place of (\ref{LP Condition for entering N}), we require:
\begin{equation} \label{QP Condition for entering N}
\begin{split}
\max_j (\Theta_{\mathcal{N}})_{jj} \| A A^T\|_{\infty} \leq \text{reg}_{thr,k}, \\
\max_j (\Theta_{\mathcal{N}})_{jj} \| Q Q^T\|_{\infty} \leq \text{reg}_{thr,k},
\end{split}
\end{equation}
\noindent where $\text{reg}_{thr,k}$ is updated as indicated in the linear programming case (sub-section \ref{subsubsection: regularization matrices LP}). As before, by permuting the columns so that the first $n_1$ correspond to indices in $\mathcal{N}$ while the remaining ones correspond to indices in $\mathcal{B}$, the augmented system in (\ref{Augmented System}) takes the form:
\begin{equation} \label{QP permuted augmented system}
\begin{bmatrix} 
-(Q_{\mathcal{N}} +\Theta_{\mathcal{N}}^{-1}+R_{p\mathcal{N}}) & -Q_{\mathcal{B}\mathcal{N}}^T & A_{\mathcal{N}}^T \\
-Q_{\mathcal{B}\mathcal{N}} &  -(Q_{\mathcal{B}} +\Theta_{\mathcal{B}}^{-1}+R_{p\mathcal{B}})& A_{\mathcal{B}}^T\\
A_{\mathcal{N}} & A_{\mathcal{B}} & R_d
\end{bmatrix}
\begin{bmatrix}
\Delta x_{\mathcal{N}}\\
\Delta x_{\mathcal{B}} \\ 
\Delta y
\end{bmatrix}
= 
\begin{bmatrix}
\xi_{d\mathcal{N}}\\
\xi_{d\mathcal{B}}\\
\xi_p
\end{bmatrix},
\end{equation}
\noindent where
\[\xi_{d\mathcal{N}} = c_{\mathcal{N}} + \begin{pmatrix} Q_{\mathcal{N}} & Q_{\mathcal{B}\mathcal{N}}^T \end{pmatrix} \begin{pmatrix}
           x_{\mathcal{N},k} \\
           x_{\mathcal{B},k} 
         \end{pmatrix} - A_{\mathcal{N}}^T y_k - \sigma_{k} \mu_k X_{\mathcal{N}}^{-1} e_{\mathcal{N}}, \]
\[\xi_{d\mathcal{B}} = c_{\mathcal{B}} + \begin{pmatrix} Q_{\mathcal{B}\mathcal{N}} & Q_{\mathcal{B}}\end{pmatrix} \begin{pmatrix}
           x_{\mathcal{N},k} \\
           x_{\mathcal{B},k} 
         \end{pmatrix} - A_{\mathcal{B}}^T y_k - \sigma_{k} \mu_k X_{\mathcal{B}}^{-1} e_{\mathcal{B}},\]
\[\xi_p = b-Ax_k,\]
\noindent and the permuted matrix $Q$ is:
\[Q = \begin{bmatrix}
Q_{\mathcal{N}} & Q_{\mathcal{B}\mathcal{N}}^T\\
Q_{\mathcal{B}\mathcal{N}} & Q_{\mathcal{B}}
\end{bmatrix},\]
\noindent with $Q_{\mathcal{N}} \in \mathbb{R}^{n_1 \times n_1}$, $Q_{\mathcal{B}\mathcal{N}} \in \mathbb{R}^{n_2 \times n_1}$ and $Q_{\mathcal{B}} \in \mathbb{R}^{n_2 \times n_2}$ being the respective blocks of the matrix $Q$, while $R_{p\mathcal{N}} \in \mathbb{R}^{n_1 \times n_1}$ and $R_{p\mathcal{B}} \in \mathbb{R}^{n_2 \times n_2}$ are the only two non-zero blocks of the block-diagonal primal regularization matrix $R_p$. As we mentioned earlier, when we solve general convex quadratic programming problems, it is dangerous to eliminate the (1,1) block of (\ref{Augmented System}) and solve the problem using (\ref{Normal Equations}), since the latter system may become dense. However in the linear programming case, our regularization matrix was tuned based on the properties of the normal equations. In order to overcome this problem, we introduce a primal regularization that can absorb the non-diagonal elements of the (1,1) block of the permuted augmented system (\ref{QP permuted augmented system}). This allows us to safely (from the sparsity and computational point of view) pivot on this block and perform the analysis in a similar manner as in the linear programming case. Hence, we define:
\begin{equation} \label{R_pN definition}
R_{p\mathcal{N}} =\big( \Delta_{p\mathcal{N}} -\text{off}(Q_{\mathcal{N}}) \big),
\end{equation}
\noindent with
\begin{equation} \label{QP Delta_pN definition}
\Delta_{p\mathcal{N}} = \|Q_{\mathcal{N}}\|_{\infty} I_{n_1},
\end{equation}
\noindent where $\Delta_{p\mathcal{N}} \in \mathbb{R}^{n_1 \times n_1}$ is a uniform diagonal matrix, which ensures that $R_{p\mathcal{N}} \succ 0$ and diagonally dominant. Although $\Delta_{p\mathcal{N}}$ can have sizeable values, (\ref{QP Condition for entering N}) ensures that the respective elements in $\Theta_{\mathcal{N}}^{-1}$ have significantly larger values, making this perturbation acceptable. Using (\ref{R_pN definition}), the (1,1) block of (\ref{QP permuted augmented system}) becomes:
\[-(Q_{\mathcal{N}} +\Theta_{\mathcal{N}}^{-1}+R_{p\mathcal{N}}) = -(\Theta_{\mathcal{N}}^{-1} + D_{p\mathcal{N}}),\]
\noindent where $D_{p\mathcal{N}} = \text{diag}(Q_{\mathcal{N}}) + \Delta_{p\mathcal{N}}$ is a diagonal matrix. For simplicity of notation, let 
\[\bar{Q}_{\mathcal{N}} = \Theta_{\mathcal{N}}^{-1} + D_{p\mathcal{N}}.\]
\noindent Pivoting on the (1,1) block of (\ref{QP permuted augmented system}) results in the following partially reduced augmented system:
\begin{equation} \label{QP partially reduced augmented system}
\begin{bmatrix} 
 Q_{\mathcal{B}\mathcal{N}}\bar{Q}_{\mathcal{N}}^{-1}Q_{\mathcal{B}\mathcal{N}}^T-(Q_{\mathcal{B}} +\Theta_{\mathcal{B}}^{-1}+R_{p\mathcal{B}})&\quad A_{\mathcal{B}}^T-Q_{\mathcal{B}\mathcal{N}}\bar{Q}_{\mathcal{N}}^{-1}A_{\mathcal{N}}^T\\
A_{\mathcal{B}}-A_{\mathcal{N}}\bar{Q}_{\mathcal{N}}^{-1}Q_{\mathcal{B}\mathcal{N}}^T & R_d + A_{\mathcal{N}}\bar{Q}_{\mathcal{N}}^{-1}A_{\mathcal{N}}^T
\end{bmatrix}
\begin{bmatrix}
\Delta x_{\mathcal{B}} \\ 
\Delta y
\end{bmatrix}
= 
\begin{bmatrix}
\xi_1\\
\xi_2
\end{bmatrix},
\end{equation}
\noindent where:
\begin{equation*}
\begin{split}
\xi_1 =&\ \xi_{d\mathcal{B}} - Q_{\mathcal{B}\mathcal{N}}\bar{Q}_{\mathcal{N}}^{-1}\xi_{d\mathcal{N}},\\
\xi_2 =&\ \xi_p + A_{\mathcal{N}}\bar{Q}_{\mathcal{N}}^{-1}\xi_{d\mathcal{N}}.
\end{split}
\end{equation*}
\noindent Using a similar reasoning as before, we will tune the matrix $R_{p\mathcal{B}}$ so that sparsity is promoted. By looking at the (1,1) block of (\ref{QP partially reduced augmented system}), one can see that an obvious choice for this matrix would be:
\begin{equation} \label{R_pB definition}
R_{p\mathcal{B}} = \big( \Delta_{p\mathcal{B}} + \text{off}(Q_{\mathcal{B}\mathcal{N}}\bar{Q}_{\mathcal{N}}^{-1}Q_{\mathcal{B}\mathcal{N}}^T) \big),
\end{equation}
\noindent with
\begin{equation} \label{QP Delta_pB definition}
\Delta_{p\mathcal{B}} = \max_j (\bar{Q}_{\mathcal{N}}^{-1})_{jj} \|Q_{\mathcal{B}\mathcal{N}} Q_{\mathcal{B}\mathcal{N}}^T \|_{\infty} I_{n_2},
\end{equation}
\noindent where $\Delta_{p\mathcal{B}} \in \mathbb{R}^{n_2 \times n_2}$ is a uniform diagonal matrix, which ensures that $R_{p\mathcal{B}} \succ 0$ and diagonally dominant. Finally, by looking at the (2,2) block of (\ref{QP partially reduced augmented system}), we can define $R_d$ in a similar manner as in the linear programming case as:
\begin{equation} \label{QP R_d definition}
R_d = \big( \Delta_d -\text{off}(A_{\mathcal{N}}\bar{Q}_{\mathcal{N}}^{-1}A_{\mathcal{N}}^T) \big),
\end{equation}
\noindent with
\begin{equation} \label{QP Delta_d definition}
\Delta_d = \max_j (\bar{Q}_{\mathcal{N}}^{-1})_{jj} \|A_{\mathcal{N}} A_{\mathcal{N}}^T \|_{\infty} I_m,
\end{equation}
\noindent where again $\Delta_d \in \mathbb{R}^{m \times m}$ is a uniform diagonal matrix, which ensures that $R_d \succ 0$ and diagonally dominant. Note that condition (\ref{QP Condition for entering N}), which defines columns qualified to enter $\mathcal{N}$, ensures that the positive elements of the diagonal matrices $\Delta_{p\mathcal{B}},\ \Delta_d$ will be strictly less than $\text{reg}_{thr,k}$, at every iteration $k$ of the algorithm. 

\paragraph*{Motivation}
$\ $\\
\noindent As in the linear programming case, let us provide the motivation for the previously presented regularization scheme. We will derive some useful bounds that extend those provided in the \textit{motivation} paragraph for the linear programming regularization. All the bounds stated here are direct applications of the results obtained in \cite{paper_9} and for simplicity are given without proofs. Let:

\begin{equation*}
M = \begin{bmatrix}
-Q-\Theta^{-1} & A^T\\
A & 0
\end{bmatrix},\ \ 
E = \begin{bmatrix}
\Delta_p & 0\\
0 &\quad \delta_d I_m
\end{bmatrix},
\end{equation*}

\noindent and denote by $\lambda_i$ and $\tilde{\lambda}_i$ the i-th smallest eigenvalues of $M$ and $M +E$ respectively.  Note that $\Delta_p$ is a permuted $n\times n$ diagonal matrix, comprised of the two uniform primal regularization matrices $\delta_{p\mathcal{N}}I_{n_1},\ \delta_{p\mathcal{B}}I_{n_2}$, with $n_1 + n_2 = n$. Let $\zeta_i = \min_{\mu \in \lambda(-Q-\Theta^{-1})} |\lambda_i - \mu|$, where $\lambda_i \in \lambda(M)$, $\lambda_i \not\in \lambda(-Q-\Theta^{-1})$ and $\lambda_i \neq 0$. Let also $Mx = \lambda_i x$, with $\|x\| =1$. Partitioning $x = [x_1^H\ x_2^H]^H$, it can be proven as before that:
\[\|x_1\| \leq \frac{\|A\|}{\sqrt{\zeta_i^2+\|A\|^2}},\ \|x_2\| \leq \frac{\|A\|}{\sqrt{\lambda_i^2+\|A\|^2}}.\]
\noindent A counterpart of Lemma \ref{lemma LP motivation x2(t)} for this case follows from \cite{paper_9} and states that if $|\lambda_i| > \delta_d + \|E\|$, then, $\forall\ t \in [0,1]$:
\[\|x_2(t)\| \leq \frac{\|A\|}{\sqrt{(|\lambda_i|  -\delta_d - \|E\|)^2 + \|A\|^2}} = \phi_i.\]
\noindent Similarly, if $\zeta_i > \|\Delta_p\| + \|E\|$, then $\forall\ t \in [0,1]$ we have:
\[\|x_1(t)\| \leq \frac{\|A\|}{\sqrt{(\zeta_i - \|\Delta_p\| - \|E\|)^2+\|A\|^2 }} = \varphi_i,\]
\noindent where $x(t) = [x_1(t)^H\ x_2(t)^H]^H$ solves the problem $(M+tE)x(t) = \lambda_i(t) x(t)$, for some $t \in [0,1]$. For a detailed derivation of the previous results, the interested reader can look at \cite{paper_9}, Lemmas 2.8, 2.9. Finally, the counterpart of Theorem \ref{Theorem LP reg motivation} for this case states that for each $i$ such that: $|\lambda_i| > \delta_d + \|E\|$ and $\zeta_i > \|\Delta_p\| + \|E\|$, we have:
\[|\lambda_i - \tilde{\lambda}_i| \leq \|\Delta_p\| \varphi_i^2 + \delta_d \phi_i^2.\]

\noindent These bounds are slightly less intuitive than the ones provided for the linear programming case, however similar arguments to those used in the linear programming case can be employed here, supporting the claim that the uniform regularization that we introduce does not perturb the sufficiently large (in the absolute value sense) eigenvalues of the non-regularized system significantly. The main reason why we provide these bounds is for completeness. We could proceed by showing, as in the linear programming case, that further dropping $\text{off}(Q_{\mathcal{N}}),\ \text{off}(Q_{\mathcal{B}\mathcal{N}}\bar{Q}_{\mathcal{N}}^{-1}Q_{\mathcal{B}\mathcal{N}}^T)$ and $\text{off}(A_{\mathcal{N}}\bar{Q}_{\mathcal{N}}^{-1}A_{\mathcal{N}}^T)$ (as the proposed non-diagonal regularization suggests) alters the eigenvalues of the diagonally regularized system in a controlled way, but for ease of presentation we omit this for a future study.

\paragraph*{Rank deficient matrices and the value of $\epsilon$}
$\ $\\
Notice that both in the linear and the quadratic programming case, during some iterations of the IPM, no columns will satisfy the respective conditions for entering $\mathcal{N}$. In order to ensure that rank deficiency will not get in the way of the proposed method, at every such iteration $k$, we apply a uniform dual regularization $R_d  = \text{reg}_{thr,k} I_m$, where $\text{reg}_{thr,k}$ is updated as stated in sub-section \ref{subsubsection: regularization matrices LP}. In the quadratic programming case, we also include a uniform primal regularization $R_p = \text{reg}_{thr,k}I_n$. We expect that sufficiently large (in the absolute value sense) eigenvalues ($\gg 2 \cdot \text{reg}_{thr,k}$) of the system are perturbed insignificantly by using such a uniform regularization. Once at least one column enters $\mathcal{N}$, we drop this uniform regularization, and start using the regularization matrices presented in this paper.\\

\noindent Notice that $\text{reg}_{thr}$ is not allowed to decrease more than a pre-specified value $\epsilon > 0$. We set this to: $\epsilon = \max \{\frac{0.1 \cdot \text{tol}}{\|A\|^2},10^{-13}\}$, where $\text{tol}$ is the error tolerance for successful termination of the algorithm and is usually set to the values $10^{-6}$ or $10^{-8}$. This value is based on the bounds derived in the motivation paragraphs presented both for the linear and the quadratic programming case, so that $\epsilon \phi_i^2$ is small.\\
$\ $\\

\section{Spectral Analysis} \label{section: Spectral analysis}

\noindent This section focuses on analysing the spectral properties of the regularized systems provided in the previous section. As before, the analysis is split into linear and quadratic programming respectively. For each of these cases, we will provide the spectral properties of the respective augmented and partially reduced augmented system, showing the effectiveness of the proposed regularization method. 

\subsection{Linear Programming}

\noindent For linear programming problems, we employ only dual regularization, that is we set $R_p = 0$ and use only $R_d \succ 0$. In sub-section \ref{subsubsection: regularization matrices LP}, it was noted that $\Delta_{d}$ is chosen such that $R_d \succ 0$ and diagonally dominant. This is very easy to see, by looking at the definition of (\ref{LP Delta_d definition}) combined with (\ref{dual regularization matrix}). Since $R_d$ is diagonally dominant, we are able to invoke the \textit{Gershgorin Circle Theorem}, which states that if:
\[r_i = \sum_{j=1,j\neq i}^{m} | (A_{\mathcal{N}} \Theta_{\mathcal{N}} A_{\mathcal{N}}^T)_{ij} |,\]
\noindent then any eigenvalue of $R_d$ is positive and lies in one of the following discs:
\[\{\lambda: |\lambda - \delta_d| \leq r_i \},\]
\noindent where $\delta_d$ is defined in (\ref{LP Delta_d definition}), $i = 1,\cdots,m$. This yields: $0 < \lambda_i \leq \delta_d + r_i, \ \forall\ i = 1,\cdots,m$, where $\lambda_i$ represents the i-th eigenvalue of $R_d$. On the other hand, by construction, we know that $\delta_d \geq r_i + \min_{j: (A_{\mathcal{N}}\Theta_{\mathcal{N}} A_{\mathcal{N}}^T)_{jj} > 0}  ( (A_{\mathcal{N}}\Theta_{\mathcal{N}} A_{\mathcal{N}}^T)_{jj})$, $\forall \ i = 1,\cdots,m$  and hence:
\begin{equation} \label{LP Gershgorin Circle Theorem bound}
\min_{j: (A_{\mathcal{N}}\Theta_{\mathcal{N}} A_{\mathcal{N}}^T)_{jj} > 0}  ( (A_{\mathcal{N}}\Theta_{\mathcal{N}} A_{\mathcal{N}}^T)_{jj}) \leq \lambda_{i} \leq \delta_d + r_i < 2 \delta_d,\ \ \ \ \forall\ i = 1,\cdots, m.
\end{equation}

\noindent Let us now analyse the spectral properties of the matrix in (\ref{LP permuted augmented system}). For that we provide the following theorem, which gives bounds for the eigenvalues of the system. The proof is based on the developments in \cite{paper_3} and \cite{paper_11}.
\begin{theorem} \label{Thm: LP Spectral Properties augmented}
For all $(x,z) > 0$ and $R_d$ as defined in \textnormal{(\ref{dual regularization matrix})}, the coefficient matrix of \textnormal{(\ref{LP permuted augmented system})} has exactly $n$ negative and $m$ positive eigenvalues. Order and denote them as:
\[\mu_{-n} \leq \mu_{-n+1} \leq \cdots \mu_{-1} < 0 < \mu_1 \leq \cdots \mu_m.\]
These eigenvalues satisfy the following bounds:
\begin{flalign*}
\mu_{-1} < -\min_j(\Theta^{-1})_{jj}, &&
\end{flalign*}
\begin{flalign*}\mu_{-n} \geq \frac{1}{2}\bigg(\big(\lambda_{\min}(R_d) - \max_j (\Theta^{-1})_{jj} \big) - \big[(\max_j(\Theta^{-1})_{jj} + \lambda_{\min}(R_d))^2 + 4(\sigma_{\max}(A))^2\big]^{\frac{1}{2}}\bigg),&&
\end{flalign*}
\begin{flalign*} \mu_{m} \leq \frac{1}{2}\bigg(2 \delta_d + \big(4 \delta_d^2 + 4(\sigma_{\max}(A))^2\big)^{\frac{1}{2}} \bigg),&&
\end{flalign*}
\begin{flalign*} \mu_{1} \geq \frac{1}{2}\bigg((\lambda_{\min}(R_d)-\max_j(\Theta^{-1})_{jj}) + \big[(\max_j(\Theta^{-1})_{jj}+\lambda_{\min}(R_d))^2 + 4(\sigma_{\min}(A))^2\big]^{\frac{1}{2}} \bigg).&&
\end{flalign*}
\noindent In case $\text{rank} (A) < m$, the eigenspace of the eigenvalues originating only from $R_d$ is $\{0\}\times \text{Null}(A^T)$ and there are $m - \text{rank}(A)$ such eigenvalues.
\end{theorem}

\noindent {\textbf{Proof}}
\noindent Firstly, from Sylvester's law of inertia we know that since $\Theta$ and $A \Theta A^T + R_d$ are positive definite, the regularized augmented system matrix of (\ref{LP permuted augmented system}) possesses precisely $n$ negative and $m$ positive eigenvalues. If $\mu$ is an eigenvalue of the linear system matrix of (\ref{LP permuted augmented system}), then there are vectors $u \in \mathbb{R}^n$ and $p \in \mathbb{R}^m$ that cannot both be zero, using which the eigenvalue problem can be written in the following form:
\begin{equation} \label{LP eigenvalue problem 1}
\begin{split}
- \Theta^{-1}u + A^T p = \mu u, \\
Au + R_d p  = \mu p.
\end{split}
\end{equation}
\noindent As observed in \cite{paper_2}, if $\text{rank}(A) < m$, there are some eigenvalues of the matrix in (\ref{LP permuted augmented system}), that satisfy: $R_d p = \mu p$. The associated eigenspace is $\{0\}\times \text{Null}(A^T)$.\\

\noindent If $\mu < 0$ then $u \neq 0$ since otherwise $p = 0$ because $R_d \succ 0$. On the other hand, if $\mu > 0$ then  $p \neq 0$ since otherwise $u = 0$ because $\Theta^{-1} \succ 0$. Taking the inner product of the first equation of (\ref{LP eigenvalue problem 1}) with $u$, and the second equation with $p$ and subtracting the former from the latter gives:

\[u^T\Theta^{-1}u + p^T R_d p = -\mu u^T u + \mu p^T p.\]

\noindent Using the fact that $\Theta^{-1} \succ 0$, along with $R_d \succ 0$, and assuming that $\mu < 0$ (i.e. $u \neq 0$):

\[(\min_j(\Theta^{-1})_{jj} + \mu) u^T u \leq \mu p^Tp,\]

\noindent  where the inequality follows because the left hand side is as small as possible and we dropped the positive term $p^T R_d p$. But since $\mu <0 $ in this case, we know that $-\min_j(\Theta^{-1})_{jj} > \mu =\mu_{-1} $. Furthermore, if $\mu < 0$ then we know that $R_d-\mu I \succ 0$. Hence it is invertible and we can solve the second equation of (\ref{LP eigenvalue problem 1}) with respect to $p$, substitute the result in the first equation and take the inner product with $u$ to get:
\[ p = - (R_d- \mu I)^{-1} Au,\]
\[-u^T \Theta^{-1} u - u^T A^T(R_d- \mu I)^{-1} A u = \mu u^T u.\]

\noindent Hence:

\[ -\max_j(\Theta^{-1})_{jj} - (\sigma_{\max}(A))^2 (\lambda_{\min}(R_d) - \mu)^{-1} \leq \mu,\]

\noindent where we observed that the left hand side has negative terms, took the most negative possible values for these terms and divided by $u^Tu$. Note that for the second term of the left hand side, we used the fact that for two positive definite matrices $A,\ B$, we have that $\lambda_{\min}(A+B) \geq \lambda_{\min}(A) + \lambda_{\min}(B)$.  Solving the previous inequality with respect to $\mu$ (and using the roots of the second order equation), we get that:

\[\mu_{-n} \geq \frac{1}{2}\bigg(\big(\lambda_{\min}(R_d) - \max_j (\Theta^{-1})_{jj} \big) - \big[(\max_j(\Theta^{-1})_{jj} + \lambda_{\min}(R_d))^2 + 4(\sigma_{\max}(A))^2\big]^{\frac{1}{2}}\bigg).\]

\noindent  Now, for the case where $\mu > 0$ (where we know that $p \neq 0$), we solve the first equation of (\ref{LP eigenvalue problem 1}) with respect to $u$, substitute the result in the second one and take the inner product with $p$, to get:

\[ u = \frac{1}{\mu}(\frac{1}{\mu} \Theta^{-1} + I)^{-1} A^T p, \]
\[ \frac{1}{\mu}p^T A(\frac{1}{\mu} \Theta^{-1} + I)^{-1} A^Tp + p^T R_d p = \mu p^T p.\]

\noindent Observe that $\lambda_{\max} ((\frac{1}{\mu} \Theta^{-1} + I)^{-1}) \leq 1$. Given that all the terms on the left hand side are positive, we can take upper bounds for every term, multiply everything by $\mu$ (since $\mu >0$) and divide both sides by $p^Tp$. This gives us the following second order inequality with respect to $\mu$:

\[\mu^2 - \lambda_{\max}(R_d)\mu - (\sigma_{\max}(A))^2 \leq 0.\]
\noindent Solving the previous quadratic inequality, gives:
\[ \mu_{m} \leq \frac{1}{2}\bigg(2 \delta_d + \big(4 \delta_d^2 + 4(\sigma_{\max}(A))^2\big)^{\frac{1}{2}} \bigg),\]

\noindent where we used the right-most upper bound given in (\ref{LP Gershgorin Circle Theorem bound}). Working similarly using the same equation but slightly altered, that is:
\[ u = ( \Theta^{-1} + \mu I)^{-1} A^T p, \]
\[ p^T A(\Theta^{-1} + \mu I)^{-1} A^Tp + p^T R_d p = \mu p^T p,\]

\noindent and by taking lower bounds on each term of the left hand side and re-arranging them, we get the following inequality:
\[\mu^2 + (\max_j(\Theta^{-1})_{jj} - \lambda_{\min}(R_d))\mu - (\sigma_{\min}(A)^2 + \max_j(\Theta^{-1})_{jj}\lambda_{\min}(R_d)) \geq 0.\]
\noindent Solving the previous gives us the last bound:
\[\mu_{1} \geq \frac{1}{2}\bigg((\lambda_{\min}(R_d)-\max_j(\Theta^{-1})_{jj}) + \big[(\max_j(\Theta^{-1})_{jj}+\lambda_{\min}(R_d))^2 + 4(\sigma_{\min}(A))^2\big]^{\frac{1}{2}} \bigg),\]
\noindent which completes the proof.
\qed

\noindent Below we provide an analogous theorem applied to the matrix of (\ref{LP partially reduced augmented system}). Again, we use the definition of $R_d$ that is given in (\ref{dual regularization matrix}). With this in mind, we know that the (2,2) block of (\ref{LP partially reduced augmented system}) is comprised of two diagonal matrices, i.e.: \[D^* = \text{diag}(A_{\mathcal{N}} \Theta_{\mathcal{N}} A_{\mathcal{N}}^T) + \Delta_d,\]
\noindent where $\Delta_d$ is defined in (\ref{LP Delta_d definition}). The proof is similar to that of the previous theorem, and hence it is not provided here.
\begin{theorem} \label{Thm: LP Spectral Properties partially reduced augmented}
For all $(x,z) > 0$ and $R_d$ as defined in \textnormal{(\ref{dual regularization matrix})}, the coefficient matrix of \textnormal{(\ref{LP partially reduced augmented system})} has exactly $n_2$ negative and $m$ positive eigenvalues. Order and denote them as:
\[\bar{\mu}_{-n_2} \leq \bar{\mu}_{-n_2+1} \leq \cdots \bar{\mu}_{-1} < 0 < \bar{\mu}_1 \leq \cdots \bar{\mu}_m.\]
These eigenvalues satisfy the following bounds:
\begin{flalign*}
\bar{\mu}_{-1} < -\min_j(\Theta_{\mathcal{B}}^{-1})_{jj},&&
\end{flalign*}
\begin{flalign*}
\bar{\mu}_{-n_2} \geq \frac{1}{2}\bigg(\big(\min_i D^*_{ii} - \max_j(\Theta_{\mathcal{B}}^{-1})_{jj} \big) - \big[(\max_j(\Theta_{\mathcal{B}}^{-1})_{jj} + (\min_i D^*_{ii}))^2 + 4(\sigma_{\max}(A_{\mathcal{B}}))^2  \big]^{\frac{1}{2}}\bigg),&&
\end{flalign*}
\begin{flalign*}
 \bar{\mu}_{m} \leq \frac{1}{2}\bigg(\max_i D^*_{ii} + \big((\max_i D^*_{ii})^2 + 4(\sigma_{\max}(A_{\mathcal{B}}))^2\big)^{\frac{1}{2}} \bigg),&&
 \end{flalign*}
\begin{flalign*}
\bar{\mu}_{1} \geq \frac{1}{2}\bigg((\min_i D^*_{ii}-\max_j(\Theta_{\mathcal{B}}^{-1})_{jj}) + \big[(\max_j(\Theta_{\mathcal{B}}^{-1})_{jj}+\min_i D^*_{ii})^2 + 4(\sigma_{\min}(A_{\mathcal{B}}))^2\big]^{\frac{1}{2}} \bigg).&&
\end{flalign*}
\noindent In case $\text{rank} (A_{\mathcal{B}}) < m$, the eigenspace of the eigenvalues originating only from $D^*$ is $\{0\}\times \text{Null}(A_{\mathcal{B}}^T)$ and there are $m - \text{rank}(A_{\mathcal{B}})$ such eigenvalues.
\end{theorem}
\noindent Now we can compare the bounds given in Theorems \ref{Thm: LP Spectral Properties augmented} and \ref{Thm: LP Spectral Properties partially reduced augmented} and observe clear advantages of using the partially reduced augmented system (\ref{LP partially reduced augmented system}) over the full augmented system (\ref{LP permuted augmented system}). Firstly, one can easily note that $-\min_j(\Theta_{\mathcal{B}}^{-1})_{jj} = -\min_j(\Theta^{-1})_{jj}$, hence the bound for the largest negative eigenvalue is identical for both systems. However, there are two main differences:
\begin{enumerate}
\item We have that $\max_j(\Theta_{\mathcal{B}}^{-1})_{jj} \leq \max_j(\Theta^{-1})_{jj}
$ (and usually $\max_j(\Theta_{\mathcal{B}}^{-1})_{jj} \ll \max_j(\Theta^{-1})_{jj}
$). As a consequence the bound on the most negative eigenvalue of (\ref{LP permuted augmented system}) will be larger (in the absolute value sense), than the bound on the respective eigenvalue of (\ref{LP partially reduced augmented system}).
\item Our guaranteed lower bound for the minimum eigenvalue of $R_d$ is smaller than the respective lower bound for the minimum eigenvalue of $D^*$. In fact,
\[ \min_i D^*_{ii} \geq \delta_d ,\]
\[ \lambda_{\min}(R_d) \geq \min_{j: (A_{\mathcal{N}}\Theta_{\mathcal{N}} A_{\mathcal{N}}^T)_{jj} > 0}  ((A_{\mathcal{N}}\Theta_{\mathcal{N}} A_{\mathcal{N}}^T)_{jj}),\]
\noindent where  $\delta_d$ is defined in (\ref{LP Delta_d definition}) and the second lower bound is given in (\ref{LP Gershgorin Circle Theorem bound}). By construction the first bound is better. As a consequence, the smallest positive eigenvalue of (\ref{LP partially reduced augmented system}) is guaranteed to be at least as large as $\delta_d$.
\end{enumerate}
\subsection{Quadratic Programming}
\noindent For quadratic programming problems we employ a primal-dual regularization. In subsection \ref{subsubsection: regularization matrices QP}, it was noted that $\Delta_{d}$ is chosen such that $R_d \succ 0$ and diagonally dominant, while $\Delta_{p\mathcal{B}}$ is chosen such that $R_{p\mathcal{B}} \succ 0$ and diagonally dominant. This can be seen by looking at (\ref{QP Delta_d definition}) combined with (\ref{QP R_d definition}) and (\ref{QP Delta_pB definition}) combined with (\ref{R_pB definition}), respectively. Similarly, positive definiteness and diagonal dominance of $R_{p\mathcal{N}}$ follows immediately by construction, i.e. by looking at equations (\ref{R_pN definition}) and (\ref{QP Delta_pN definition}). For notational convenience, we define:
\[\bar{Q}_{\mathcal{N}} = \Theta_{\mathcal{N}}^{-1} + \text{diag}(Q_{\mathcal{N}}) + \Delta_{p\mathcal{N}}.\]
\begin{itemize}
\item For $R_d$, we are able to invoke the Gershgorin circle theorem as in the linear programming case stating that if:
\[r_i = \sum_{j=1,j\neq i}^{m} | (A_{\mathcal{N}} \bar{Q}_{\mathcal{N}}^{-1} A_{\mathcal{N}}^T)_{ij} |,\]
\noindent then any eigenvalue of $R_d$ is positive and lies in one of the following discs:
\[\{\lambda: |\lambda - \delta_d| \leq r_i \},\]
\noindent where $\delta_d = \max_{j}(\bar{Q}_{\mathcal{N}}^{-1})_{jj} \|A_{\mathcal{N}} A_{\mathcal{N}}^T\|_{\infty} $, $i = 1,\cdots,m$. This yields: $0 < \lambda_{i} \leq \delta_d + r_i, \ \forall\ i = 1,\cdots,m$, where $\lambda_i$ is the i-th eigenvalue of $R_d$. On the other hand, by construction we know that $\delta_d \geq r_i + \min_{j: (A_{\mathcal{N}}\bar{Q}_{\mathcal{N}}^{-1} A_{\mathcal{N}}^T)_{jj} > 0}  ((A_{\mathcal{N}}\bar{Q}_{\mathcal{N}}^{-1} A_{\mathcal{N}}^T)_{jj})$, $\forall \ i = 1,\cdots,m$ and hence:
\begin{equation} \label{QP Gershgorin Circle Theorem bound R_d}
\min_{j:(A_{\mathcal{N}}\bar{Q}_{\mathcal{N}}^{-1} A_{\mathcal{N}}^T)_{jj} > 0}  ((A_{\mathcal{N}}\bar{Q}_{\mathcal{N}}^{-1} A_{\mathcal{N}}^T)_{jj}) \leq \lambda_{i} \leq \delta_d + r_i < 2 \delta_d.
\end{equation}
\item  For $R_{p\mathcal{B}}$, we apply the same theorem however in this case we have:
\[r_i = \sum_{j=1,j\neq i}^{n_2} | (Q_{\mathcal{B}\mathcal{N}} \bar{Q}_{\mathcal{N}}^{-1} Q_{\mathcal{B}\mathcal{N}}^T)_{ij} |,\]
\noindent and any eigenvalue of $R_{p\mathcal{B}}$ is positive and lies in one of the following discs:
\[\{\lambda: |\lambda - \delta_{pB}| \leq r_i \},\]
\noindent where $\delta_{p\mathcal{B}} = \max_{j}(\bar{Q}_{\mathcal{N}}^{-1})_{jj} \|Q_{\mathcal{B}\mathcal{N}} Q_{\mathcal{B}\mathcal{N}}^T\|_{\infty} $, $i = 1,\cdots,n_2$.  As before, we know that:
\begin{equation} \label{QP Gershgorin Circle Theorem bound R_pB}
\min_{j:(Q_{\mathcal{B}\mathcal{N}}\bar{Q}_{\mathcal{N}}^{-1} Q_{\mathcal{B}\mathcal{N}}^T)_{jj} > 0}  ((Q_{\mathcal{B}\mathcal{N}}\bar{Q}_{\mathcal{N}}^{-1}Q_{\mathcal{B}\mathcal{N}}^T)_{jj}) \leq \lambda_{i} \leq \delta_{p\mathcal{B}} + r_i < 2 \delta_{p\mathcal{B}},\ \  \forall\ i =1,\cdots,n_2=|\mathcal{B}|,
\end{equation}
\noindent where $\lambda_i$ is the i-th eigenvalue of $R_{p\mathcal{B}}$.
\item  Finally, we can work similarly to examine the spectral properties of $R_{p\mathcal{N}}$. Again by letting:
\[r_i = \sum_{j=1,j\neq i}^{n_1} (Q_{\mathcal{N}})_{ij},\]
\noindent any eigenvalue of $R_{p\mathcal{N}}$ is positive and lies in one of the following discs:
\[\{\lambda:|\lambda-\delta_{p\mathcal{N}}| \leq r_i \},\]
\noindent  where $\delta_{p\mathcal{N}} = \|Q_{\mathcal{N}}\|_{\infty}$, $i = 1,\cdots,n_1$. This yields: $0 < \lambda_i  \leq \delta_{p\mathcal{N}} + r_i$, $\forall\ i =1,\cdots,n_1$, where $\lambda_i$ is the i-th eigenvalue of $R_{p\mathcal{N}}$. But since $Q_{\mathcal{N}} \succeq 0$ as a principal minor of $Q \succeq 0$, we know that if a diagonal element of $Q_{\mathcal{N}}$ is zero, then its respective column and row are also zero. Hence this implies tighter final bounds, that is:
\begin{equation} \label{QP Gershgorin Circle Theorem bound R_pN}
\min_{j:(Q_{\mathcal{N}})_{jj} > 0} ( (Q_{\mathcal{N}})_{jj}) < \lambda_i  \leq \delta_{p\mathcal{N}} + r_i < 2\delta_{p\mathcal{N}},\ \  \forall\ i =1,\cdots,n_1=|\mathcal{N}|.
\end{equation}
\end{itemize}
\noindent Let us now analyse the spectral properties of (\ref{QP permuted augmented system}). For that we provide the following theorem, which is the extension of Theorem \ref{Thm: LP Spectral Properties augmented} for the QP case. The proof is almost identical and hence it is not provided here. For notational convenience, let: 
\[H = Q + \Theta^{-1} + R_p.\]
\begin{theorem} \label{Thm: QP Spectral Properties augmented}
For all $(x,z) > 0$ and $R_d,\ R_{p\mathcal{B}},\ R_{p\mathcal{N}}$ as defined in \textnormal{(\ref{QP R_d definition}), (\ref{R_pB definition})} and \textnormal{(\ref{R_pN definition})} respectively, the coefficient matrix of \textnormal{(\ref{QP permuted augmented system})} has exactly $n$ negative and $m$ positive eigenvalues. Order and denote them as:
\[\mu_{-n} \leq \mu_{-n+1} \leq \cdots \mu_{-1} < 0 < \mu_1 \leq \cdots \mu_m.\]
These eigenvalues satisfy the following bounds:
\begin{flalign*}
\mu_{-1} < -\lambda_{\min}(H), &&
\end{flalign*}
\begin{flalign*}
\mu_{-n} \geq \frac{1}{2}\bigg(\big(\lambda_{\min}(R_d) - \lambda_{\max}(H) \big) - \big[(\lambda_{\max}(H) + \lambda_{\min}(R_d))^2 + 4(\sigma_{\max}(A))^2 \big]^{\frac{1}{2}}\bigg),&&
\end{flalign*}
\begin{flalign*} \mu_{m} \leq \frac{1}{2}\bigg(2 \delta_d + \big(4 \delta_d^2 + 4(\sigma_{\max}(A))^2\big)^{\frac{1}{2}} \bigg),&&
\end{flalign*}
\begin{flalign*} \mu_{1} \geq \frac{1}{2}\bigg((\lambda_{\min}(R_d)-\lambda_{\max}(H)) + \big[(\lambda_{\max}(H)+\lambda_{\min}(R_d))^2 + 4(\sigma_{\min}(A))^2\big]^{\frac{1}{2}} \bigg).&&
\end{flalign*}
\noindent In case $\text{rank} (A) < m$, the eigenspace of the eigenvalues originating only from $R_d$ is $\{0\}\times \text{Null}(A^T)$ and there are $m - \text{rank}(A)$ such eigenvalues.
\end{theorem}

\noindent Below we provide a similar theorem, applied to (\ref{QP partially reduced augmented system}). For that, we will use $R_d,\ R_{p\mathcal{B}},\ R_{p\mathcal{N}}$ as defined in sub-section \ref{subsubsection: regularization matrices QP} as well as the respective eigenvalue bounds given in (\ref{QP Gershgorin Circle Theorem bound R_d}), (\ref{QP Gershgorin Circle Theorem bound R_pB}) and (\ref{QP Gershgorin Circle Theorem bound R_pN}). Using the definitions of the regularization matrices, we know that the matrix in the (1,1) block of (\ref{QP partially reduced augmented system}) takes the form:
\[\bar{H} = (Q_{\mathcal{B}} + \Theta_{\mathcal{B}}^{-1} + \Delta_{p\mathcal{B}} - \text{diag}(Q_{\mathcal{B}\mathcal{N}}\bar{Q}_{\mathcal{N}}^{-1}Q_{\mathcal{B}\mathcal{N}}^T)),\]
\noindent while the (2,2) block of (\ref{QP partially reduced augmented system}) becomes:
\[D^* = \text{diag}(A_{\mathcal{N}} \bar{Q}_{\mathcal{N}}^{-1} A_{\mathcal{N}}^T) + \Delta_d.\]
\begin{theorem} \label{Thm: QP Spectral Properties partially reduced augmented}
For all $(x,z) > 0$ and $R_d,\ R_{p\mathcal{B}},\ R_{p\mathcal{N}}$ as defined in \textnormal{(\ref{QP R_d definition}), (\ref{R_pB definition})} and \textnormal{(\ref{R_pN definition})} respectively, the coefficient matrix of \textnormal{(\ref{QP partially reduced augmented system})} has exactly $n_2$ negative and $m$ positive eigenvalues. Order and denote them as:
\[\bar{\mu}_{-n_2} \leq \bar{\mu}_{-n+1} \leq \cdots \bar{\mu}_{-1} < 0 < \bar{\mu}_1 \leq \cdots \bar{\mu}_m.\]
These eigenvalues satisfy the following bounds:
\begin{flalign*}
\bar{\mu}_{-1} < -\lambda_{\min}(\bar{H}),  &&
\end{flalign*}
\begin{flalign*} 
\bar{\mu}_{-n_2} \geq &\ \frac{1}{2}\bigg(\big(\min_{j}(D^*)_{jj} - \lambda_{\max}(\bar{H}) \big)  -  \\ &\ -\big[(\lambda_{\max}(\bar{H}) + \min_{j}(D^*)_{jj})^2 + 4(\sigma_{\max}(A_{\mathcal{B}}-A_{\mathcal{N}}\bar{Q}_{\mathcal{N}}^{-1}Q_{\mathcal{B}\mathcal{N}}^T))^2\big]^{\frac{1}{2}}\bigg), &&
\end{flalign*}
\begin{flalign*}
\bar{\mu}_{m} \leq \frac{1}{2}\bigg(\max_{j}(D^*)_{jj} + \big(\max_{j}(D^*)_{jj}^2 + 4(\sigma_{\max}(A_{\mathcal{B}}-A_{\mathcal{N}}\bar{Q}_{\mathcal{N}}^{-1}Q_{\mathcal{B}\mathcal{N}}^T))^2\big)^{\frac{1}{2}} \bigg),&& 
\end{flalign*}
 \begin{flalign*}
 \bar{\mu}_{1} \geq &\ \frac{1}{2}\bigg(( \min_{j}(D^*)_{jj}-\lambda_{\max}(\bar{H}))\ + \\ &\ + \big[(\lambda_{\max}(\bar{H})+ \min_{j}(D^*)_{jj})^2 + 4(\sigma_{\min}(A_{\mathcal{B}}-A_{\mathcal{N}}\bar{Q}_{\mathcal{N}}^{-1}Q_{\mathcal{B}\mathcal{N}}^T))^2\big]^{\frac{1}{2}} \bigg).&&
 \end{flalign*}
 \noindent In case $\text{rank} (A_{\mathcal{B}}-A_{\mathcal{N}}\bar{Q}_{\mathcal{N}}^{-1}Q_{\mathcal{B}\mathcal{N}}^T) < m$, the eigenspace of the eigenvalues originating only from $D^*$ is $\{0\} \times \text{Null} (A_{\mathcal{B}}^T-Q_{\mathcal{B}\mathcal{N}}\bar{Q}_{\mathcal{N}}^{-1}A_{\mathcal{N}}^T)$ and there are $m - \text{rank} (A_{\mathcal{B}}-A_{\mathcal{N}}\bar{Q}_{\mathcal{N}}^{-1}Q_{\mathcal{B}\mathcal{N}}^T)$ such eigenvalues.
\end{theorem}
\noindent Let us compare the bounds given in Theorems \ref{Thm: QP Spectral Properties augmented} and \ref{Thm: QP Spectral Properties partially reduced augmented} to observe once again the advantages of using the partially reduced augmented system (\ref{QP partially reduced augmented system}) over the full augmented system (\ref{QP permuted augmented system}). There are three significant differences in the eigenvalue bounds of these two systems:
\begin{enumerate}
\item For the bound on the largest negative eigenvalue of the two systems, we know that:
 \[\lambda_{\min}(H) \geq \min_j (\Theta^{-1})_{jj} + \lambda_{\min}(R_p),\]
 \noindent where
 \[ \lambda_{\min}(R_p) \geq \min\bigg\{ \min_{j: (Q_{\mathcal{B}\mathcal{N}}\bar{Q}_{\mathcal{N}}^{-1} Q_{\mathcal{B}\mathcal{N}}^T)_{jj} > 0}  ( (Q_{\mathcal{B}\mathcal{N}}\bar{Q}_{\mathcal{N}}^{-1}Q_{\mathcal{B}\mathcal{N}}^T)_{jj}),\ \min_{j: (Q_{\mathcal{N}})_{jj} > 0} ( (Q_{\mathcal{N}})_{jj}) \bigg\},\]
 \noindent from (\ref{QP Gershgorin Circle Theorem bound R_pB}) and (\ref{QP Gershgorin Circle Theorem bound R_pN}) respectively. However, since $\min_j (\Theta_{\mathcal{B}}^{-1})_{jj} \ll \min_j (\Theta_{\mathcal{N}}^{-1})_{jj}$ we can conclude that:
 \[\lambda_{\min}(H) \geq \min_j (\Theta_{\mathcal{B}}^{-1})_{jj} + \lambda_{\min}(R_p),\]
 \noindent while
 \[\lambda_{\min}(\bar{H}) \geq \min_j (\Theta_{\mathcal{B}}^{-1})_{jj} +\max_{j}(\bar{Q}_{\mathcal{N}}^{-1})_{jj} \|Q_{\mathcal{B}\mathcal{N}} Q_{\mathcal{B}\mathcal{N}}^T\|_{\infty} - \max_j (Q_{\mathcal{B}\mathcal{N}} \bar{Q}_{\mathcal{N}}^{-1} Q_{\mathcal{B}\mathcal{N}}^T)_{jj},\]
\noindent where we used (\ref{QP Delta_pB definition}) as the definition of $\Delta_{p\mathcal{B}}$. We observe that the difference: \[\max_{j}(\bar{Q}_{\mathcal{N}}^{-1})_{jj} \|Q_{\mathcal{B}\mathcal{N}} Q_{\mathcal{B}\mathcal{N}}^T\|_{\infty} - \max_j (Q_{\mathcal{B}\mathcal{N}} \bar{Q}_{\mathcal{N}}^{-1} Q_{\mathcal{B}\mathcal{N}}^T)_{jj},\]
\noindent increases as more elements enter the set $\mathcal{N}$. On the other hand, $\lambda_{\min}(R_p)$ is expected to decrease at every iteration of the interior-point method. Hence the bound on $\bar{\mu}_{-1}$ is expected to be better than that on $\mu_{-1}$, as more elements enter the partition $\mathcal{N}$.

\item For the bound on the most negative eigenvalue of the two systems, we know that:
\[\lambda_{\max}(H) \leq \lambda_{\max}(Q) +\max_j (\Theta^{-1})_{jj} + \lambda_{\max}(R_p),\]
 \noindent where $\lambda_{\max}(R_p) \leq 2\max\{\delta_{p\mathcal{N}},\delta_{p\mathcal{B}}\}$. However, since $\max_j (\Theta_{\mathcal{N}}^{-1})_{jj} \geq \max_j (\Theta_{\mathcal{B}}^{-1})_{jj}$, we observe that: 
 \[\lambda_{\max}(H) \leq \lambda_{\max}(Q) + \max_j (\Theta_{\mathcal{N}}^{-1})_{jj} + \lambda_{\max}(R_p),\]
 \noindent where we used the definition of $\Delta_{p\mathcal{N}}$ given in (\ref{QP Delta_pN definition}). On the other hand, 
 \[\lambda_{\max}(\bar{H})\leq \lambda_{\max}(Q_{\mathcal{B}}) +  \max_j (\Theta_{\mathcal{B}}^{-1})_{jj} + (\Delta_{p\mathcal{B}})_{ii},\quad \forall\ i \in \{1,\cdots,n\}.\]
 \noindent where, from (\ref{QP Delta_pB definition}), we know that $(\Delta_{p\mathcal{B}})_{ii} = \max_j (\bar{Q}_{\mathcal{N}}^{-1})_{jj} \|Q_{\mathcal{B}\mathcal{N}} Q_{\mathcal{B}\mathcal{N}}^T\|_{\infty},\ \forall\ i \in \{1,\cdots,n\}$. Clearly the bound on $\lambda_{\max}(\bar{H})$ is significantly smaller than that on $\lambda_{\max}(H)$, since it is usually the case that $\max_j (\Theta_{\mathcal{N}}^{-1})_{jj} \gg \max_j (\Theta_{\mathcal{B}}^{-1})_{jj}$, while $\lambda_{\max}(R_p) > \max_j (\bar{Q}_{\mathcal{N}}^{-1})_{jj} \|Q_{\mathcal{B}\mathcal{N}} Q_{\mathcal{B}\mathcal{N}}^T\|_{\infty}$. Hence, the most negative eigenvalue of (\ref{QP partially reduced augmented system}) is expected to have a significantly smaller magnitude than that of (\ref{QP permuted augmented system}).
 
 \item As in the LP case, our guaranteed lower bound for the minimum eigenvalue of $R_d$ is smaller than the respective lower bound for the minimum eigenvalue of $D^*$. In fact,
\[ \min_i D^*_{ii} \geq \delta_d\]
\[ \lambda_{\min}(R_d) \geq \min_{j: (A_{\mathcal{N}}\bar{Q}_{\mathcal{N}}^{-1} A_{\mathcal{N}}^T)_{jj} > 0}  ((A_{\mathcal{N}}\bar{Q}_{\mathcal{N}}^{-1} A_{\mathcal{N}}^T)_{jj}),\]
\noindent where we use $\delta_d$ as defined in (\ref{QP Delta_d definition}), while the last inequality follows from (\ref{QP Gershgorin Circle Theorem bound R_d}). By construction, the first bound is better. As a consequence, the smallest positive eigenvalue of (\ref{QP partially reduced augmented system}) is guaranteed to be at least as large as $\delta_d$. 
\end{enumerate}

\section{Implementation and Numerical Results} \label{section: Implementation}
\subsection{The Algorithmic Framework}

\noindent At this point, we are providing a generic algorithm (\ref{Algorithm IPM}), summarizing the infeasible primal-dual IPM with non-diagonal regularization. The algorithm solves the Newton system arising from the optimality conditions of (\ref{Primal Barrier Problem})-(\ref{Dual Barrier Problem}), at each iteration, using a direct method. Note that this is just a general outline and does not contain the actual details of the implemented method. Implementation details will be presented in the next sub-section. Note that in the algorithm, we make the distinction between linear and quadratic programming problems, by using the logical variables LP and QP, respectively.\\

\begin{singlespace}
\begin{algorithm} 
\caption{Infeasible primal-dual IPM with non-diagonal regularization}
    \label{Algorithm IPM}
    \textbf{Input:}  $A, Q, b, c$, $\text{tol}$, $\text{maxit}$\\
    \textbf{Parameters:} $0< \sigma_{\min} \leq \sigma_{\max}$ (bounds for the centring parameter), $\epsilon = \max\big\{\frac{\text{tol}\cdot 10^{-1}}{\|A\|_2^2},10^{-13} \big\}.$ \\
    \textbf{Initial point:} Choose a well-centred $w_0 = (x_0,y_0,r_0,s_0,z_0)$  with $x_0,z_0 \geq 0$, $\mu_0 = \frac{x_0^T z_0}{n}$, $k = 0$, $\text{reg}_{thr,0} = 1$.
\begin{algorithmic}
\State $\text{res}_p^0 = b - Ax_0,\ \text{res}_d^0 = c - A^Ty_0 - z_0 + Qx_0 $. 
\While {($k < \text{maxit}$)}
\If {(($\|\text{res}_p^k\| < \text{tol}) \wedge (\|\text{res}_d^k\| < \text{tol}) \wedge (\mu_k < \text{tol})$)}
	\State Declare convergence and return the optimal solution.
	\State \Return $(x_k,y_k,z_k)$.	
\Else
	\State  $\text{reg}_{thr,k} = \max \{ O(\mu_k), \epsilon \}$.
	\If {($\mathcal{N} = \emptyset$)} 
		\State $R_d = \text{reg}_{thr,k} I_m$.
		\If {(QP)}
			\State $R_p =  \text{reg}_{thr,k} I_n$.
		\EndIf
	\Else
		\If {(LP)}
			\State $R_d$ from (\ref{dual regularization matrix}) and (\ref{LP Delta_d definition}), $R_p = 0$.
		\ElsIf {(QP)}
			\State $R_d$ from (\ref{QP R_d definition}), (\ref{QP Delta_d definition}) and $R_p$ from (\ref{R_pN definition}), (\ref{QP Delta_pN definition}), (\ref{R_pB definition}), (\ref{QP Delta_pB definition}).
		\EndIf 
	\EndIf \\
	
    \State Choose $\sigma_k \in [\sigma_{\min},\sigma_{\max}].$
	\If {(LP)}
		\State Compute $\Delta w_k = (\Delta x_k,\Delta y_k,\Delta r_k,\Delta z_k)$ by solving (\ref{Normal Equations}) and the substitutions.
		\State ($s_k = 0, \Delta s = 0$).
	\ElsIf {(QP)}    
        \State Compute $\Delta w_k = (\Delta x_k, \Delta y_k,\Delta r_k,\Delta s_k,\Delta z_k)$ by solving (\ref{QP partially reduced augmented system}) and the substitutions.
    \EndIf
    \State $a_x^{\max} = \text{min}_{\Delta x_i < 0} \bigg \{1,-\frac{x_i}{\Delta x_i}\bigg\},\ \ a_z^{\max} = \text{min}_{\Delta z_i < 0} \bigg \{1,-\frac{z_i}{\Delta z_i}\bigg\}.$
    \State $x_k(a) = x_k + \tau a_x^{\max} \Delta x, r_k(a) = r_k + \tau a_x^{\max}\Delta r$.
    \State $z_k(a) = z_k + \tau a_z^{\max} \Delta z,\ y_k(a) = y_k + \tau a_z^{\max} \Delta y,\ s_k(a) = s_k + \tau a_z^{\max} \Delta s$, $\ \ \ \tau \in\ ]0,1[$. 
    \State $\mu_k(a) = \frac{x_k(a)^Tz_k(a)}{n}$.
    \State $k = k+1$.
\EndIf
\EndWhile 
\end{algorithmic}
\end{algorithm}
\end{singlespace}
\subsection{Implementation Details}

\noindent We implemented the algorithm in Matlab. Our implementation solves linear and convex quadratic programming problems in the standard form. However, all the free variables are treated as variables bounded by some box constraints. We set some initial bounds, $$l_f = -10^{2} \leq x_f \leq 10^2 = u_f,$$
\noindent for all the free variables. If the method pushes some of these variables to take values outside of this box, then the respective bounds are increased to give space for variables to increase their values. Note that this heuristic causes that extra iterations are needed to converge for a few problems, since every time the box constraints are changed, the method loses primal feasibility.

\paragraph*{Regularization} 
$\ $\\
\noindent We set $\text{reg}_{thr,0} = 1$, and we decrease it at the same rate as $\mu_k$ decreases, until it becomes smaller than $\epsilon = \max\big\{\frac{\text{tol}\cdot 10^{-1}}{\|A\|_2^2},10^{-13} \big\}$. Then, it takes this value and stays constant for the rest of the optimization process. As before, $\text{tol}$ is the error tolerance specified by the user. At every iteration, we enable columns to enter the set $\mathcal{N}$ only if: $\max_{j \in \mathcal{N}} (\Theta)_{jj} \max\big\{\|AA^T\|_{\infty},\|QQ^T\|_{\infty}\big\} \leq \text{reg}_{thr,k}$. This ensures that $(\Delta_{d})_{ii}$, as defined in (\ref{LP Delta_d definition}) and (\ref{QP Delta_d definition}) for linear and convex quadratic problems respectively, is smaller than $\text{reg}_{thr,k}, \forall i \in \{1,\cdots,m\},\ \forall\ k \geq 0$. The latter also holds for $(\Delta_{p\mathcal{B}})_{ii}$ as in (\ref{QP Delta_pB definition}) $, \forall i \in \{1,\cdots,n_2\}$, which is only defined for quadratic programming problems. Of course for linear programming problems we have $R_p = 0$. Note that during the first iterations of the method, $\mathcal{N}$ is usually empty. In order to avoid instability, we include a uniform dual regularization $R_d = \text{reg}_{thr,k} I_m$. For the quadratic programming case, we also include a uniform primal regularization, that is: $R_p = \text{reg}_{thr,k} I_n$. This uniform regularization is dropped when $\mathcal{N}$ is non-empty. As an extra safeguard, when the factorization of the system fails, we increase $\text{reg}_{thr}$ by a factor of 10 and repeat the factorization. If this process is repeated for 6 consecutive times, we stop the method. All other implementation details concerning the regularization follow from Section \ref{section: exact reg}.\\

\paragraph*{Newton-step computation}
$\ $\\
\noindent  For general convex quadratic problems, the Newton direction is calculated from system (\ref{QP partially reduced augmented system}), after computing its symmetric $LDL^T$ decomposition, where $L$ is a lower triangular matrix and $D$ is diagonal. For that, we use the build-in Matlab symmetric decomposition (i.e. \texttt{ldl}). We know that such a decomposition always exists, with $D$ diagonal, for the aforementioned system, since after introducing the regularization, the matrix of (\ref{QP partially reduced augmented system}) is guaranteed to be quasi-definite; a class of matrices known to be strongly factorizable, \cite{paper_4}. For that reason, we change the default pivot threshold of \texttt{ldl} to $10^{-14}$. We use such a small pivot threshold in order to avoid any 2x2 pivots during the factorization routine. For linear programming problems, we solve the system (\ref{Normal Equations}) (with $Q = 0$), using the build-in Cholesky decomposition of Matlab (i.e. \texttt{chol}). $\Delta x$ is then recovered from (\ref{Delta x recovery}). In the quadratic programming case, $\Delta s$ is recovered from (\ref{Delta s recovery}). In both cases $\Delta z$ is recovered from (\ref{Delta z recovery}) and $\Delta r$ from (\ref{Delta r recovery}).

\paragraph*{Starting point} 
$\ $\\
\noindent We have already mentioned that the method is infeasible and hence the starting point does not need to be primal and dual feasible. The only requirement is that the initial values of the variables $x,\ z$ are strictly positive. We use a starting point that was proposed in \cite{paper_6}. Here we will only state it for completeness. To construct this point, we try to solve the pair of problems (\ref{non-regularized primal}), (\ref{non-regularized dual}), but we ignore the non-negativity constraints. Such relaxed problems have closed form solutions:

\begin{equation} \label{eq: mehr start point 1}
\tilde{x} = A^T(AA^T)^{-1}b,\ \ \tilde{y} = (AA^T)^{-1}A(c+Q\tilde{x}), \ \ \tilde{z} = c - A^T \tilde{y} + Q\tilde{x}.
\end{equation}

\noindent Then, in order to guarantee positivity and sufficient magnitude of $x,z$, we compute the expressions $\delta_x = \text{max}(-1.5 \text{min}\{\tilde{x}_i\},0)$ and $\delta_z = \text{max}(-1.5 \text{min}\{\tilde{z}_i\},0)$ and we obtain:

\begin{equation}
\tilde{\delta_x} = \delta_x + 0.5\frac{(\tilde{x}+ \delta_x e)^T(\tilde{z}+ \delta_z e)}{\sum_{i=1}^n (\tilde{z}_i + \delta_z)},
\end{equation}

\begin{equation}
\tilde{\delta_z} = \delta_z + 0.5\frac{(\tilde{x}+ \delta_x e)^T(\tilde{z}+ \delta_z e)}{\sum_{i=1}^n (\tilde{x}_i + \delta_x)},
\end{equation}
\noindent where $e$ is the vector of ones of appropriate dimension. Finally, we define the starting point by setting:

\begin{equation}
r^0 = 0,\ \ s^0 = 0,\ \ y^0 = \tilde{y},\ \ z_i^0 = \tilde{z}_i + \tilde{\delta_z},\ \ x_i^0 = \tilde{x}_i + \tilde{\delta_x},\ \ \ \ i = 1 ,\cdots ,n.
\end{equation}

\paragraph*{Centring parameter}
$\ $\\
\noindent As minimum and maximum centring parameters, we fix $\sigma_{\min} = 0.05$ and $\sigma_{\max} = 0.95$. In the first iteration we use $\sigma_0 = 0.5$. Then, at each iteration $k$, in order to determine the centring parameter $\sigma_k$, we perform the following operations:
\[\sigma_k = \max \{ (1-a_x^{k-1})^5,(1-a_z^{k-1})^5 \},\]
\noindent where $a_x^{k-1}, a_z^{k-1}$ are the step-lengths in directions $\Delta x,\ \Delta z$ of the previous iteration, respectively. Then we assign:
\[\sigma_k = \min \{ \sigma_k, \sigma_{\max} \},\]
\noindent and finally 
\[\sigma_k = \max \{ \sigma_k, \sigma_{\min} \}.\]
\noindent The latter is a heuristic which performs well in infeasible IPM implementations.

\paragraph*{Step-length computation}
$\ $\\
\noindent In order to calculate the step-length, we apply the \textit{fraction to the boundary} rule, that is we compute the largest step-lengths to the boundary of the non-negative orthant, i.e.:
\begin{equation} \label{eq: Implem eq 1}
a_x^{\max} = \text{min}_{\Delta x_i < 0} \bigg \{1,-\frac{x_i}{\Delta x_i}\bigg\},\ \ a_z^{\max} = \text{min}_{\Delta z_i < 0} \bigg \{1,-\frac{z_i}{\Delta z_i}\bigg\},
\end{equation}
\noindent and we set:
\begin{equation}\label{eq: Implem eq 2}
a_x = \tau a_x^{\max},\ \ a_z= \tau a_z^{\max},
\end{equation}
\noindent where $\tau \in\ ]0,1[$ is set to $\tau =  0.995$. The constant $\tau$ acts as a safeguard against bad directions. Taking a full step towards a direction can potentially push the iterates of the algorithm close to the boundary. This in turn, can significantly slow down the convergence of the method. The primal variables $x, r$ are updated using the step-length $a_x$ while the dual variables $y,\ s,\ z$ are updated using the step-length $a_z$.

\paragraph*{Termination Criteria}
$\ $\\
\noindent Finally, the algorithm is terminated either if the number of maximum iterations specified by the user is reached, or when all the following three conditions are satisfied:
\[\frac{\| c -A^T y + Qx - z \|}{\|c\| + 1} \leq \text{tol},\]
\[ \frac{\|b - Ax\|}{\|b\| + 1} \leq \text{tol},\]
\noindent and
\[ \mu \leq \text{tol},\]
\noindent where $\text{tol}$ is the tolerance specified by the user.

\subsection{Numerical Results}

\noindent We have made a particular effort to keep the implementation as simple as possible, so that the regularization effects can easily be seen and analysed. For that reason, we applied scaling only to problems which required it to converge and this was needed only for 5 out of the 218 problems solved. On the other hand, no predictor-corrector technique was included. We tested our method on problems coming from the Netlib collection \cite{paper_7} as well as on a set of convex quadratic programming problems given in \cite{paper_8}. We present the numerical results, firstly for linear programming problems and then for quadratic programming ones. In order to demonstrate the effects of the proposed regularization method, we will compare it with an interior point method that uses a uniform regularization. This uniform regularization scheme, can be interpreted as the application of a standard proximal point method, in contrast to the proposed method, which can be interpreted as the application of a generalized proximal point method. The experiments were conducted on a PC with a 2.2GHz Inter Core i5 processor (dual-core) and 4GB RAM, run under Linux operating system. The Matlab version used was R2018a.

\paragraph*{Linear programming problems}
$\ $\\
\noindent As we have already stated, for linear programming problems we use only dual regularization, that is we set $R_p = 0$ and $s = 0$ in (\ref{Primal Problem})-(\ref{Dual Problem}). For that reason, we will compare our method with an algorithm that uses a uniform dual regularization, $R_d = \text{reg}_{thr,k} I_m,\ \forall k \geq 0$, where $\text{reg}_{thr,k}$ is updated as indicated in the previously presented \textit{Regularization} paragraph. If $\mathcal{N} = \emptyset$, the two methods are exactly the same. Hence, the difference between the methods arises when some columns of the constraint matrix have entered the set $\mathcal{N}$. The tolerance used in the experiments for the linear programming problems was $\text{tol}= 10^{-6}$. We will not use a smaller tolerance because our method does not have primal regularization. As a consequence, if some elements of $\Theta_{\mathcal{B}}$ become very large, this can create numerical instability if there is no primal regularization to keep such entries manageable in terms of machine precision. As an extra safeguard, when the factorization fails, we increase the uniform regularization value by a factor of 10 until the factorization is completed successfully. Finally, we set the maximum iterations of the method to be $\text{maxit} = 200$. If this number is reached, the algorithm stops indicating that the optimal solution was not found. To conclude we use:
\[\text{tol} = 10^{-6},\ \text{maxit} = 200.\]
\noindent The statistics of runs of the proposed IPM with non-diagonal regularization and of the previously mentioned IPM with uniform regularization, over the Netlib test set, have been collected in Table \ref{Netlib table}. Notice that Table \ref{Netlib table} contains only a sub-set of the 96 problems of the Netlib collection. All problems for which the set $\mathcal{N}$ stayed empty throughout the whole optimization process have been excluded. In this case, the two methods are completely equivalent. 

\begin{singlespace}
\centering
\begin{scriptsize}

\begin{longtable}{|p{2cm}| p{1cm}| p{2.5cm}| p{1.2cm}| p{1cm}| p{2.5cm}| p{1.2cm}|} 
\caption{\small Netlib Colletion} \label{Netlib table}\\
\endfirsthead

\multicolumn{7}{l}{\small{\textbf{\tablename \ \thetable:}\  \textit{Continued}}} \\
\hline
\textbf{Name} & \multicolumn{3}{c|}{\textbf{Non-diagonal Reg.}} & \multicolumn{3}{c|}{\textbf{Uniform Reg.}}\\ 
 \hline
  & \textbf{Iter.} & \textbf{$t_{\text{total}}$} (sec.) &\textbf{Stat.}  & \textbf{Iter.} & \textbf{$t_{\text{total}}$} (sec.) & \textbf{Stat.}\\

\hline
\endhead
\hline \multicolumn{7}{c}{} \\
\endfoot
\hline
\endlastfoot
\hline
\textbf{Name} & \multicolumn{3}{c|}{\textbf{Non-diagonal Reg.}} & \multicolumn{3}{c|}{\textbf{Uniform Reg.}}\\ 
 \hline
  & \textbf{Iter.} & \textbf{$t_{\text{total}}$} (sec.) &\textbf{Stat.}  & \textbf{Iter.} & \textbf{$t_{\text{total}}$} (sec.) & \textbf{Stat.}\\
\hline
ADLITTLE & 23 & 5.189100e-02 & opt  & 23 & 2.186500e-02  & opt \\
AFIRO & 10 & 1.149000e-02  & opt  & 10 & 7.599000e-03 & opt \\
AGG & 33 & 6.976600e-02  & opt  & 31 & 8.973400e-02  & opt  \\
AGG2 & 35 & 1.187540e-01  & opt  & 35 & 1.185300e-01 & opt  \\
AGG3 & 31 & 1.046800e-01  & opt  & 31 & 1.106120e-01  & opt \\
BEACONFD & 13 & 8.038000e-03  & opt  & 13 & 9.213000e-03 & opt\\
BNL1 & 43 & 2.571860e-01  & opt  & 43 & 2.555120e-01 & opt\\
CAPRI & 29 & 1.000590e-01 & opt  & 28 & 1.170850e-01 & opt \\
CZPROB & 43 & 2.985060e-01  & opt  & 48 & 3.398990e-01  & opt  \\
D2Q06C & 54 & 1.984321e+00  & opt  & 54 & 2.050860e+00  & opt  \\
DEGEN2 & 21 & 1.527010e-01  & opt  & 21 & 1.503810e-01 & opt \\
DFL001 & 84 & 1.069954e+01  & opt & 82 & 1.278976e+01  & opt   \\
FFFFF800 & 49 & 2.383090e-01  & opt  & 49 & 2.245300e-01  & opt  \\
FINNIS & 32 & 9.045000e-02  & opt  & 32 & 6.704100e-02 & opt \\
FIT2D & 42 & 2.229528e+00  & opt  & 42 & 2.303429e+00  & opt \\
FORPLAN & 31 & 1.041280e-01 & opt  & 31 & 1.285160e-01 & opt  \\
GANGES & 26 & 7.910500e-02  & opt  & 26 & 7.961900e-02  & opt \\
GFRD-PNC & 37 & 5.083300e-02  & opt  & 37 & 6.297500e-02 & opt  \\
GREENBEA & 69 & 1.467105e+00 & opt  & 69 & 1.715615e+00  & opt  \\
GREENBEB & 68 & 1.467883e+00  & opt  & 68 & 1.487645e+00  & opt  \\
GROW15 & 21 & 9.174400e-02  & opt  & 21 & 9.294100e-02 & opt \\
GROW22 & 22 & 1.397460e-01 & opt & 22 & 1.345540e-01  & opt  \\
GROW7 & 20 & 4.463700e-02  & opt  & 20 & 4.709500e-02  & opt \\
MAROS & 34 & 1.894710e-01 & opt  & 34 & 1.933730e-01  & opt  \\
MODSZK1 & 30 & 1.236590e-01  & opt  & 30 & 1.177030e-01  & opt   \\
NESM & 52 & 6.727080e-01 & opt   & 53 & 7.009060e-01  & opt  \\
PEROLD & 78 & 1.023070e+00  & opt  & 83 & 1.346218e+00 & opt\\
PILOT.JA & 98 & 2.484631e+00  & opt  & 149 & 5.730954e+00 & opt  \\
PILOT.WE & 84 & 5.628290e-01  & opt  & 80 & 6.386240e-01  & opt\\
QAP12 & 34 & 5.324705e+00  & opt  & 32 & 5.794378e+00  & opt  \\
QAP15 & 37 & 3.105651e+01  & opt  & 39 & 4.094642e+01  & opt \\
QAP8 & 19 & 2.675480e-01  & opt   & 20 & 3.163920e-01  & opt \\
SC50A & 12 & 5.564000e-03 & opt  & 12 & 7.275000e-03  & opt \\
SCAGR25 & 30 & 4.775800e-02  & opt  & 30 & 4.343100e-02  & opt  \\
SCAGR7 & 31 & 2.389600e-02  & opt  & 34 & 4.341000e-02 & opt \\
SCORPION & 32 & 2.650000e-02  & opt  & 36 & 3.096300e-02  & opt \\
SCSD1 & 19 & 2.114700e-02  & opt  & 21 & 4.742600e-02  & opt  \\
SCSD6 & 68 & 9.546700e-02  & opt  & 73 & 8.867700e-02  & opt \\
SCSD8 & 35 & 7.679500e-02  & opt    & 83 & 1.762570e-01  & opt\\
SEBA & 17 & 7.285000e-03  & opt   & 19 & 9.740000e-03 & opt\\
SHELL & 40 & 1.353860e-01  & opt   & 40 & 1.891540e-01 & opt  \\
SHIP04L & 26 & 1.208970e-01  & opt  & 32 & 2.571370e-01  & opt \\
SHIP04S & 30 & 9.137400e-02  & opt  & 26 & 9.265600e-02  & opt  \\
SHIP08L & 31 & 2.447510e-01 & opt  & 33 & 2.590480e-01  & opt\\
SHIP08S & 32 & 1.126960e-01  & opt & 30 & 1.044440e-01  & opt   \\
SHIP12L & 34 & 4.511480e-01  & opt & 34 & 3.488350e-01  & opt \\
SHIP12S & 33 & 1.340110e-01  & opt  & 33 & 1.306230e-01  & opt\\
SIERRA & 32 & 3.101600e-01 & opt  & 33 & 3.564710e-01 & opt  \\
VTP.BASE & 26 & 1.852400e-02 & opt   & 26 & 1.572000e-02  & opt  \\
WOOD1P & 34 & 1.002611e+00 & opt  & 33 & 1.406915e+00  & opt \\

\hline
\end{longtable}
\end{scriptsize}
\end{singlespace}
\noindent Both IPM with non-diagonal regularization and IPM with uniform regularization, solved all 96 problems of the Netlib collection. The former did so in 146,6 seconds and a total of 3322 IPM iterations. The latter needed 165,7 seconds and a total of 3442 iterations. In other words, the IPM using the proposed regularization, solved the whole set in 11.5\% less time, requiring 3\% less iterations. The computational benefits of the non-diagonal regularization become obvious in the larger instances of the Netlib collection. See for example problems DFL001, QAP15 in Table \ref{Netlib table}.\\

\noindent  We also include Table \ref{Sparsity table netlib}, in which the factorization times are compared when using non-diagonal and uniform regularization respectively, over the last four iterations of problems DFL001 and GREENBEA. The size of the respective constraint matrices also includes columns which were added to transform the problems to the standard form. Extra information, concerning the cardinality of the partition $\mathcal{N}$, the iteration count as well as the time needed to compute the Cholesky factorization of the system matrix at the respective iteration, has been collected in Table \ref{Sparsity table netlib}.

\begin{singlespace}
\centering
\begin{scriptsize}
\begin{longtable}{|p{2cm}| p{1cm}| p{1cm}| p{1cm}| p{1cm}| p{2cm}| p{1cm}| p{2cm}|}
\caption{Sparsity introduced from the non-diagonal regularization (linear programming)} \label{Sparsity table netlib}\\
\endfirsthead

\multicolumn{8}{l}{\small{\textbf{\tablename \ \thetable:}\  \textit{Continued}}} \\
\hline
\multirow{2}{*}{\textbf{Name}} & \multirow{2}{*}{$m$} &\multirow{2}{*}{$n$} & \multicolumn{3}{c|}{\textbf{Non-diagonal Reg.}} & \multicolumn{2}{c|}{\textbf{Uniform Reg.}}\\  \cline{4-8}

  & & &\textbf{Iter.} & $|\mathcal{N}|$ & \textbf{$t_{\text{fact}}$} (sec.) & \textbf{Iter.} & \textbf{$t_{\text{fact}}$} (sec.)\\
\hline
\endhead
\hline \multicolumn{8}{c}{} \\
\endfoot
\hline
\endlastfoot

\hline \multicolumn{8}{c}{} \\
\endfoot
\hline
\multirow{2}{*}{\textbf{Name}} & \multirow{2}{*}{$m$} &\multirow{2}{*}{$n$} & \multicolumn{3}{c|}{\textbf{Non-diagonal Reg.}} & \multicolumn{2}{c|}{\textbf{Uniform Reg.}}\\  \cline{4-8}

  & & &\textbf{Iter.} & $|\mathcal{N}|$ & \textbf{$t_{\text{fact}}$} (sec.) & \textbf{Iter.} & \textbf{$t_{\text{fact}}$} (sec.)\\
\hline
\multirow{4}{*}{DFL001}  & \multirow{4}{*}{9785} & \multirow{4}{*}{15477}  & 81 & 4089 & 0.0508 & 79 & 0.0952 \\* \cline{4-8}
 & & & 82 & 5709 & 0.0295 & 80 & 0.0972 \\*  \cline{4-8}
 & & & 83 &  6247 & 0.0258 & 81 & 0.0979\\* \cline{4-8}
& & & 84 & 7280 & 0.0166 & 82 & 0.0977\\* 
\hline
\multirow{4}{*}{GREENBEA}  & \multirow{4}{*}{3770} & \multirow{4}{*}{5973}  & 66 & 2512 & 0.0033 & 66 & 0.0107 \\* \cline{4-8}
 & & & 67 & 2536 & 0.0029 & 67 & 0.0106 \\*  \cline{4-8}
 & & & 68 &  1210 & 0.0080 & 68 & 0.0111\\* \cline{4-8}
& & & 69 & 2647 & 0.0026 & 69 & 0.0113 \\*
\hline
\end{longtable}
\end{scriptsize}
\end{singlespace}
\noindent Analysing the results reported in Tables \ref{Netlib table} and \ref{Sparsity table netlib}, one can observe that while the proposed non-diagonal regularization matrix does not affect the convergence of the method, it can accelerate the factorization significantly through the sparsity that it introduces in the system matrix. Notice that for both DFL001 and GREENBEA, almost half of their columns lie in the partition $\mathcal{N}$ and this does not prevent the algorithm from converging.\\

\begin{figure}[H]
\centering
\begin{subfigure}[H]{0.47\textwidth}
	\includegraphics[width=\textwidth]{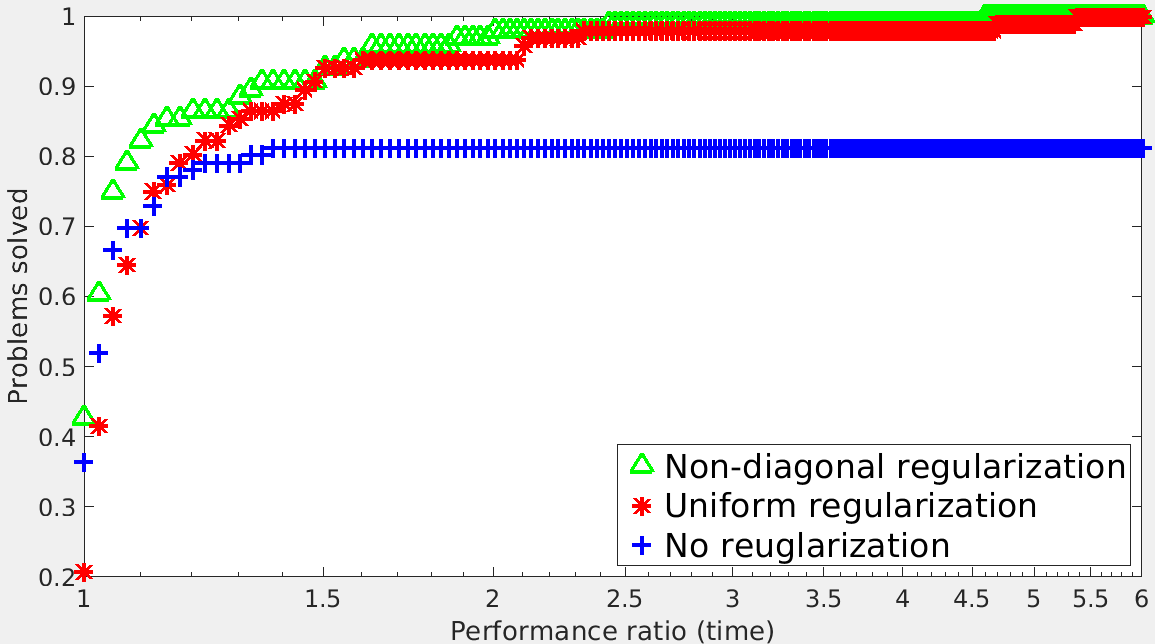}
	\caption{Performance profile with respect to time}
	\label{LP perf. prof. time}
\end{subfigure}
\quad
\begin{subfigure}[H]{0.47\textwidth}
	\includegraphics[width =\textwidth]{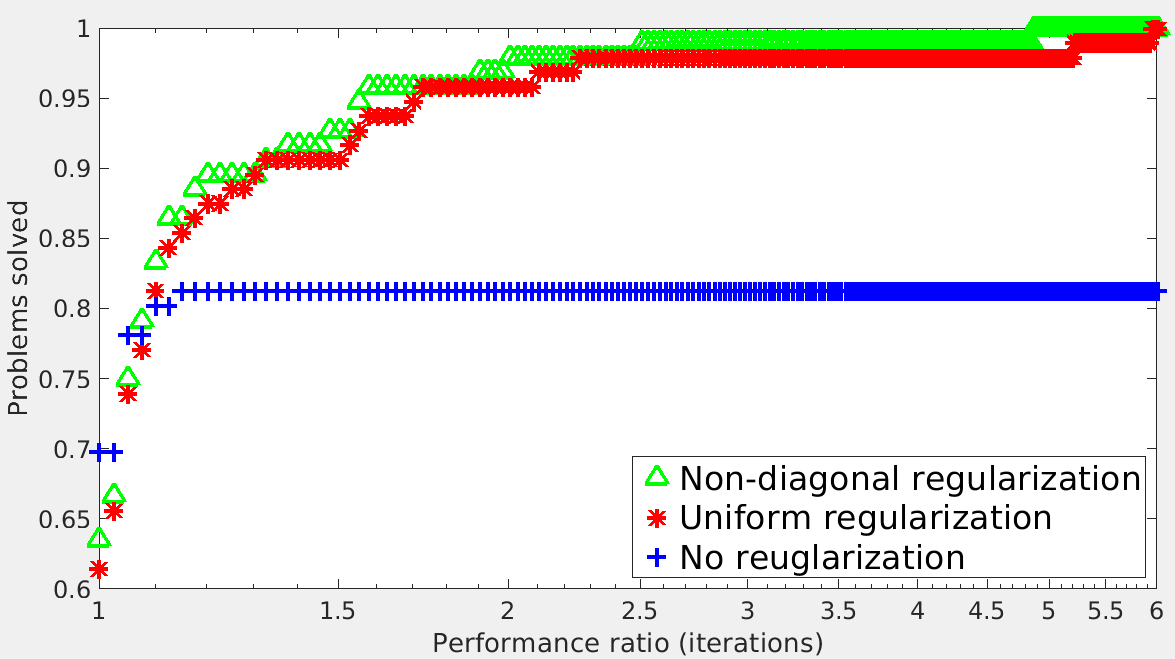}
	\caption{Performance profile with respect to iterations}
	\label{LP perf. prof. iter}
\end{subfigure}
\caption{Performance profiles over the Netlib test set.} \label{LP perf. prof.}
\end{figure}

\noindent Finally, in order to present the importance of regularization, as well as the overall comparison of the two different regularization schemes, we also include Figure \ref{LP perf. prof.}, which contains the performance profiles, over the whole Netlib set, of three different methods. The green triangles correspond to the IPM with non-diagonal regularization. The red stars correspond to the IPM with uniform regularization, and finally the blue crosses correspond to an IPM without regularization. In Figure \ref{LP perf. prof. time}, we present the performance profiles with respect to the total time to convergence, while in Figure \ref{LP perf. prof. iter} the performance profiles with respect to the total number of iterations. The horizontal axis (in logarithmic scale), represents the performance ratio with respect to the best performance achieved by one of the three methods for each problem. For example, 2 in the horizontal axis is interpreted as: ``what percentage of problems was solved by each method, in at most 2 times the best achieved time for each problem". The vertical axis shows the percentage of problems solved by each method for different values of the performance ratio. Efficiency is measured by the rate at which each of the lines increases, as the ratio increases. Robustness is measured by the maximum percentage achieved by each of the methods. For more information about performance profiles, we refer the reader to \cite{paper_27}, where this benchmarking scheme was proposed. \\

\noindent By looking at Figure \ref{LP perf. prof.}, one can observe the importance of regularization in terms of robustness of the method. The IPM scheme that does not employ any regularization, fails to solve 18.75\% of the problems in the Netlib collection. On the other hand, the IPM with non-diagonal regularization is more efficient in terms of time to convergence, when compared to the other two methods. Notice that this is not the case for the IPM using uniform regularization, which is less efficient than the other two methods for 70\% of the problems. As expected, the IPM that does not use regularization converges in less iterations for most of the problems that it successfully solves. This is expected, since in the regularized schemes, we are perturbing the Newton system. Obviously, this perturbation is benign, in the sense that it allows us to significantly improve the robustness of the method. 
\paragraph*{Convex quadratic programming problems}
$\ $\\
\noindent For this class of problems, we employ a primal-dual dynamic regularization. Hence, we will compare our method with an algorithm that uses a uniform primal-dual regularization. Such a method adds two uniform diagonal matrices $R_p = \text{reg}_{thr,k} I$ and $R_d = \text{reg}_{thr,k} I$ to the (1,1) and (2,2) blocks of the augmented system respectively. This scheme can be interpreted as the primal and dual application of the standard proximal point method, in contrast to the proposed regularization scheme, which is the primal and dual application of a generalized proximal point method. As an extra safeguard, when the factorization fails, we increase the uniform regularization value by a factor of 10 until the factorization is completed successfully. The tolerance used in the experiments for this class of problems was $\text{tol}= 10^{-8}$. As in the linear programming case, we set the maximum iterations of the method to be $\text{maxit} = 200$. To conclude we use:
\[\text{tol} = 10^{-8},\ \text{maxit} = 200.\]
\noindent For this problem set, the algorithm did not employ any scaling in the problem matrices. The computational results, obtained with the proposed non-diagonally regularized IPM and with the previously mentioned uniformly regularized IPM, over the Maros and M\'esz\'aros repository of convex quadratic programming problems, are presented in Table \ref{Maros table}. As before, Table \ref{Maros table} contains only a sub-set of the 122 problems of the collection. All problems for which the set $\mathcal{N}$ stayed empty throughout the whole optimization process have been excluded. 

\begin{singlespace}
\centering
\begin{scriptsize}
\begin{longtable}{|p{2cm}| p{1cm}| p{2.5cm}| p{1.2cm}| p{1cm}| p{2.5cm}| p{1.2cm}|}
\caption{Maros-M\'esz\'aros repository of convex quadratic problems} \label{Maros table}\\
\endfirsthead

\multicolumn{7}{l}{\small{\textbf{\tablename \ \thetable:}\  \textit{Continued}}} \\
\hline
\textbf{Name} & \multicolumn{3}{c|}{\textbf{Non-diagonal Reg.}} & \multicolumn{3}{c|}{\textbf{Uniform Reg.}}\\ 
 \hline
  & \textbf{Iter.} & \textbf{$t_{\text{total}}$} (sec.) &\textbf{Stat.}  & \textbf{Iter.} & \textbf{$t_{\text{total}}$} (sec.) & \textbf{Stat.}\\

\hline
\endhead
\hline \multicolumn{7}{c}{} \\
\endfoot
\hline
\endlastfoot
\hline
\textbf{Name} & \multicolumn{3}{c|}{\textbf{Non-diagonal Reg.}} & \multicolumn{3}{c|}{\textbf{Uniform Reg.}}\\ 
 \hline
  & \textbf{Iter.} & \textbf{$t_{\text{total}}$} (sec.) &\textbf{Stat.}  & \textbf{Iter.} & \textbf{$t_{\text{total}}$} (sec.) & \textbf{Stat.}\\
\hline
AUG2D      & 23 & 6.061128e+00  &opt & 23 & 5.894982e+00 & opt \\
AUG2DCQP   & 20 & 1.709764e+00  &opt & 20 & 1.552568e+00 & opt \\
AUG2DQP    & 23 & 1.790057e+00  &opt & 23 & 1.789244e+00  & opt \\
AUG3D & 24 & 5.630410e-01 & opt & 24 & 5.587120e-01 & opt \\
AUG3DCQP & 14 & 1.478520e-01  &opt & 14 & 1.444330e-01  &opt\\
AUG3DQP & 19 & 2.092870e-01  &opt & 19 & 1.823560e-01  &opt \\
CVXQP1     & 57 & 1.021041e+02  &opt & 57 & 1.052619e+02  &opt \\
CVXQP2     & 36 & 3.016804e+01 &opt & 36 & 3.127249e+01  &opt \\
CVXQP2     & 26 & 2.486000e-02 & opt  & 27 & 3.391610e-01  & opt \\
CVXQP3     & 38 & 1.289052e+02  &opt & 38 & 1.335597e+02  &opt \\
CVXQP3     & 36 & 6.324500e-01  &opt & 36 & 6.770390e-01  &opt \\
DTOC3      & 51 & 4.833583e+00  &opt & 51 & 4.825231e+00  &opt\\
GENHS28    & 24 & 1.257300e-02  &opt & 24 & 1.168300e-02  &opt\\
GOULDQP3   & 18 & 9.434200e-02  &opt & 32 & 1.012600e-01  &opt \\
HS118      & 21 & 1.024400e-02  &opt  & 21 & 1.006200e-02  &opt \\
HS268 & 33 & 2.154400e-02  &opt & 33 & 2.258500e-02  &opt \\
HUES-MOD   & 41 & 2.463455e+00  &opt & 41 & 2.511943e+00  &opt \\
HUESTIS    & 47 & 2.757462e+00  &opt & 46 & 3.387538e+00  &opt\\
KSIP       & 19 & 1.412779e+00  &opt  & 19 & 2.117894e+00  &opt \\
LISWET1    & 23 & 2.675218e+00 &opt & 23 & 2.810544e+00  &opt \\
LISWET10   & 48 & 4.722692e+00  &opt & 48 & 5.033838e+00  &opt \\
LISWET11   & 41 & 4.221731e+00 &opt & 41 & 4.392309e+00  &opt \\
LISWET12   & 79 & 7.336671e+00  &opt & 79 & 7.665449e+00  &opt\\
LISWET2    & 26 & 2.961113e+00  &opt & 26 & 3.077866e+00  &opt \\
LISWET5    & 38 & 3.931599e+00 &opt  & 38 & 4.055212e+00 &opt\\
LISWET7    & 34 & 3.599977e+00  &opt & 34 & 3.757427e+00  &opt\\
LISWET8    & 88 & 8.331742e+00  &opt & 86 & 8.234124e+00 &opt\\
LISWET9    & 92 & 8.660267e+00 &opt & 92 & 8.667601e+00 &opt\\
LOTSCHD    & 16 & 7.799000e-03  &opt & 16 & 5.648000e-03  &opt \\
MOSARQP1   & 18 & 1.203620e-01  &opt & 18 & 1.127110e-01  &opt \\
MOSARQP2   & 19 & 9.990200e-02 &opt & 19 & 9.517100e-02  &opt \\
POWELL20   & 34 & 3.942514e+00  &opt & 34 & 4.004035e+00  &opt\\
25FV47 & 50 & 5.708470e+00  &opt & 50 & 6.315691e+00  &opt \\
ADLITTLE & 26 & 2.163400e-02 &opt & 26 & 2.271800e-02  &opt \\
AFIRO & 23 & 1.025900e-02  &opt & 23 & 1.471700e-02  &opt\\
BANDM & 34 & 9.784200e-02  &opt & 34 & 1.229650e-01  &opt \\
BEACONFD & 26 & 6.674200e-02 &opt & 26 & 9.206800e-02  &opt \\
BORE3D & 35 & 1.082910e-01  &opt & 35 & 7.955600e-02 &opt\\
BRANDY & 32 & 1.072950e-01 &opt & 32 & 7.607900e-02  &opt \\
CAPRI & 83 & 2.571410e-01  &opt & 82 & 2.757540e-01  &opt \\
ETAMACRO & 49 & 4.680490e-01  &opt & 49 & 4.650360e-01 &opt\\
FFFFF800 & 55 & 4.667720e-01  &opt & 55 & 4.816020e-01 &opt \\
FORPLAN   & 65 & 2.766850e-01  &opt & 65 & 2.764550e-01 &opt \\
GFRD-PNC & 50 & 1.540650e-01  &opt & 50 & 1.561170e-01  &opt \\
ISRAEL  & 46 & 1.318150e-01  &opt  & 46 & 1.452930e-01  &opt\\
QPCBLEND   & 31 & 2.477900e-02  &opt & 33 & 2.679200e-02 &opt\\
QPCBOEI1   & 37 & 1.778260e-01  &opt & 37 & 1.845410e-01 &opt \\
QPCBOEI2   & 38 & 6.798500e-02  &opt & 38 & 7.008600e-02  &opt \\
QPCSTAIR   & 50 & 2.324790e-01  &opt & 50 & 2.253200e-01  &opt \\
SC205     & 25 & 2.619200e-02 &opt & 25 & 2.825100e-02  &opt \\
SCAGR25 & 34 & 6.203900e-02  &opt & 34 & 6.464000e-02 &opt \\
SCAGR7 & 32 & 2.657700e-02  &opt & 32 & 4.613000e-02  &opt\\
SCFXM1 & 40 & 1.220160e-01  &opt & 40 & 1.408360e-01  &opt \\
SCFXM2 & 49 & 2.765340e-01  &opt & 49 & 2.945260e-01  &opt \\
SCFXM3 & 50 & 4.010390e-01  &opt & 50 & 4.176580e-01  &opt \\
SCORPION & 45 & 8.277000e-02  &opt & 47 & 9.129200e-02  &opt\\
SCRS8 & 47 & 1.590580e-01  &opt & 47 & 1.773890e-01 &opt\\
SCSD1 & 23 & 5.059900e-02  &opt & 22 & 4.881200e-02 &opt \\
SCSD6 & 74 & 2.744840e-01  &opt & 67 & 2.572550e-01  &opt \\ 
SCSD8 & 21 & 1.417630e-01  &opt & 21 & 1.454910e-01  &opt \\
SCTAP1 & 32 & 5.994600e-02  &opt & 32 & 6.404600e-02 &opt \\
SEBA & 40 & 3.843170e-01  &opt & 40 & 4.018590e-01  &opt\\
SHARE2B & 41 & 4.095400e-02  &opt & 41 & 4.198200e-02  &opt \\
SHELL & 55 & 2.496066e+00  &opt & 55 & 2.906110e+00 &opt \\
SHIP04L & 25 & 1.382760e-01  &opt & 25 & 1.331710e-01  &opt \\
SHIP04S & 25 & 9.612000e-02 &opt & 25 & 9.150700e-02  &opt\\
SHIP08L & 28 & 1.215146e+00  &opt & 28 & 1.374338e+00  &opt \\
SHIP08S & 27 & 3.556720e-01  &opt & 27 & 3.986100e-01  &opt \\
SHIP12L & 32 & 2.166390e+00  &opt & 33 & 2.766354e+00  &opt\\
SHIP12S & 34 & 5.151170e-01 &opt & 34 & 5.692160e-01  &opt \\
STAIR & 32 & 1.749250e-01  &opt & 32 & 2.151530e-01  &opt \\
STANDATA & 27 & 1.789160e-01  &opt & 27 & 1.660020e-01  &opt \\ 
STCQP1     & 23 & 2.111884e+00  &opt & 23 & 2.761944e+00 &opt\\
STCQP2     & 24 & 2.103470e+00  &opt & 24 & 2.445700e+00 &opt\\
\hline
\end{longtable}
\end{scriptsize}
\end{singlespace}

\noindent In contrast to the linear programming case, the results collected in Table \ref{Maros table} do not demonstrate any significant advantage in terms of sparsity of linear systems achievable by the new regularization technique. This is a consequence of the fact that the problems under consideration are of small to medium size, while the overhead of setting up the partially reduced augmented system (\ref{QP partially reduced augmented system}) is time consuming in Matlab, where manipulating a permuted matrix is costly, due to Matlab's default mechanism to store matrices by columns. Nevertheless, both IPM with non-diagonal regularization and IPM with uniform regularization, solved all 122 problems. The former required 386,1 seconds and a total of 4162 IPM iterations. The latter required 400,2 seconds and a total of 4170 iterations. In other words, the non-diagonal scheme required 3\% less time and a similar number of iterations, as compared to the uniform scheme, for this test set.\\

\noindent As before, in order to illustrate the effect of the non-diagonal regularization in terms of factorization performance, we provide Table \ref{Sparsity table maros}, in which the factorization times obtained when using non-diagonal and uniform regularization respectively are compared, over the last four iterations of problems LISWET1, FORPLAN and SHELL. The size of the constraint matrix in each case also includes columns which were added to transform the problem to the standard form.  Information concerning the cardinality of the partition $\mathcal{N}$, the iteration count as well as the time needed to compute the $LDL^T$ factorization of the system matrix at the respective iteration, is gathered in Table \ref{Sparsity table maros}.

\begin{singlespace}
\centering
\begin{scriptsize}
\begin{longtable}{|p{2cm}| p{1cm}| p{1cm}| p{1cm}| p{1cm}| p{2cm}| p{1cm}| p{2cm}|}
\caption{Sparsity introduced from the non-diagonal regularization (quadratic programming)} \label{Sparsity table maros}\\
\endfirsthead

\multicolumn{8}{l}{\small{\textbf{\tablename \ \thetable:}\  \textit{Continued}}} \\
\hline
\multirow{2}{*}{\textbf{Name}} & \multirow{2}{*}{$m$} &\multirow{2}{*}{$n$} & \multicolumn{3}{c|}{\textbf{Non-diagonal Reg.}} & \multicolumn{2}{c|}{\textbf{Uniform Reg.}}\\  \cline{4-8}
\nopagebreak
  & & &\textbf{Iter.} & $|\mathcal{N}|$ & \textbf{$t_{\text{fact}}$} (sec.) & \textbf{Iter.} & \textbf{$t_{\text{fact}}$} (sec.)\\
  \pagebreak
\hline
\endhead
\hline \multicolumn{8}{c}{} \\
\endfoot
\hline
\endlastfoot

\hline
\multirow{2}{*}{\textbf{Name}} & \multirow{2}{*}{$m$} &\multirow{2}{*}{$n$} & \multicolumn{3}{c|}{\textbf{Non-diagonal Reg.}} & \multicolumn{2}{c|}{\textbf{Uniform Reg.}}\\  \cline{4-8}

  & & &\textbf{Iter.} & $|\mathcal{N}|$ & \textbf{$t_{\text{fact}}$} (sec.) & \textbf{Iter.} & \textbf{$t_{\text{fact}}$} (sec.)\\
\hline
\multirow{4}{*}{LISWET1}  & \multirow{4}{*}{20002} & \multirow{4}{*}{30004}  & 20 & 9670 & 0.0574 & 19 & 0.0747 \\ \cline{4-8}

 & & & 21 & 9815 & 0.0601 & 21 & 0.0692 \\  \cline{4-8}
 & & & 22 &  9935 & 0.0632 & 22 & 0.0787\\ \cline{4-8}
& & & 23 & 9984 & 0.0593 & 23 & 0.0715\\ 
\hline
\multirow{4}{*}{FORPLAN}  & \multirow{4}{*}{186} & \multirow{4}{*}{517}  & 62 & 199 & 0.0013 & 62 & 0.0036 \\* \cline{4-8}
 & & & 63 & 199 & 0.0018 & 63 & 0.0034 \\*  \cline{4-8}
 & & & 64 &  199 & 0.0013 & 64 & 0.0034 \\* \cline{4-8}
& & & 65 & 199 & 0.0016 & 65 & 0.0033 \\* 
\hline
\multirow{4}{*}{SHELL}  & \multirow{4}{*}{903} & \multirow{4}{*}{2144}  & 52 & 563 & 0.0035 & 52 & 0.0112 \\ \cline{4-8}
 & & & 53 & 565 & 0.0034 & 53 & 0.0121 \\  \cline{4-8}
 & & & 54 &  565 & 0.0033 & 54 & 0.0109 \\ \cline{4-8}
& & & 55 & 721 & 0.0033 & 55 & 0.0119 \\ 
\hline
\end{longtable}
\end{scriptsize}
\end{singlespace}

\noindent The examples presented in Table \ref{Sparsity table maros}, confirm the previous observations drawn from the linear programming examples. In particular, we can observe the benefits of the proposed non-diagonal regularization, in terms of factorization performance. On the other hand, the convergence of the method does not seem to be affected when big part of the columns of the constraint matrix lie in partition $\mathcal{N}$. \\

\begin{figure}[H]
\centering
\begin{subfigure}[H]{0.47\textwidth}
	\includegraphics[width=\textwidth]{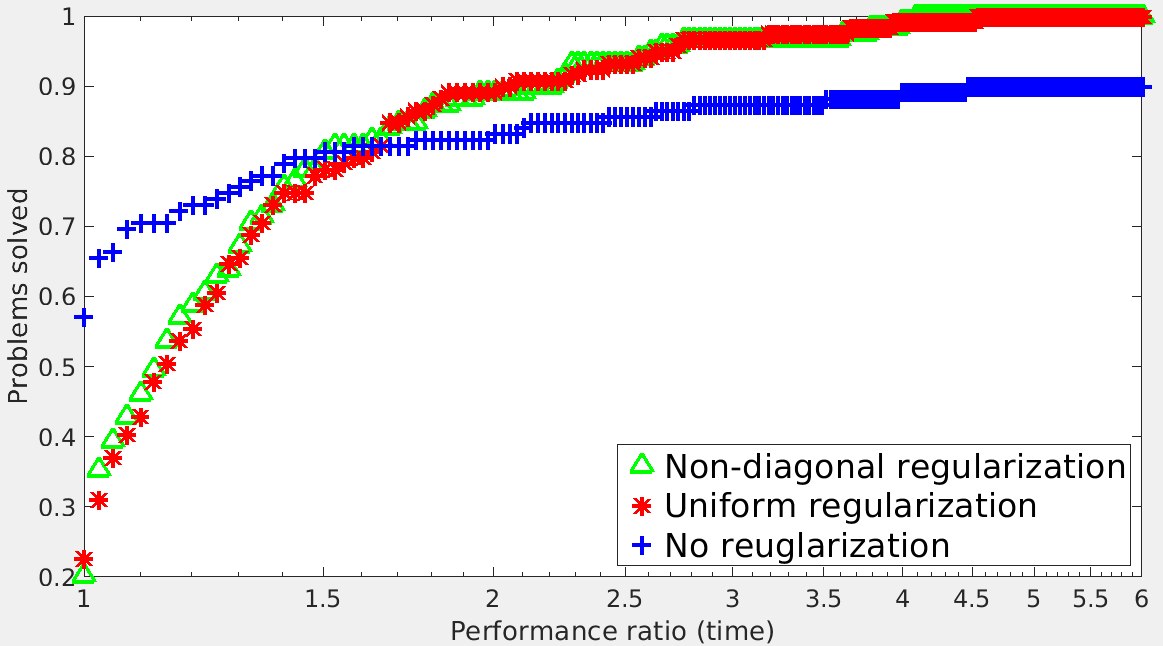}
	\caption{Performance profile with respect to time}
	\label{QP perf. prof. time}
\end{subfigure}
\quad
\begin{subfigure}[H]{0.47\textwidth}
	\includegraphics[width =\textwidth]{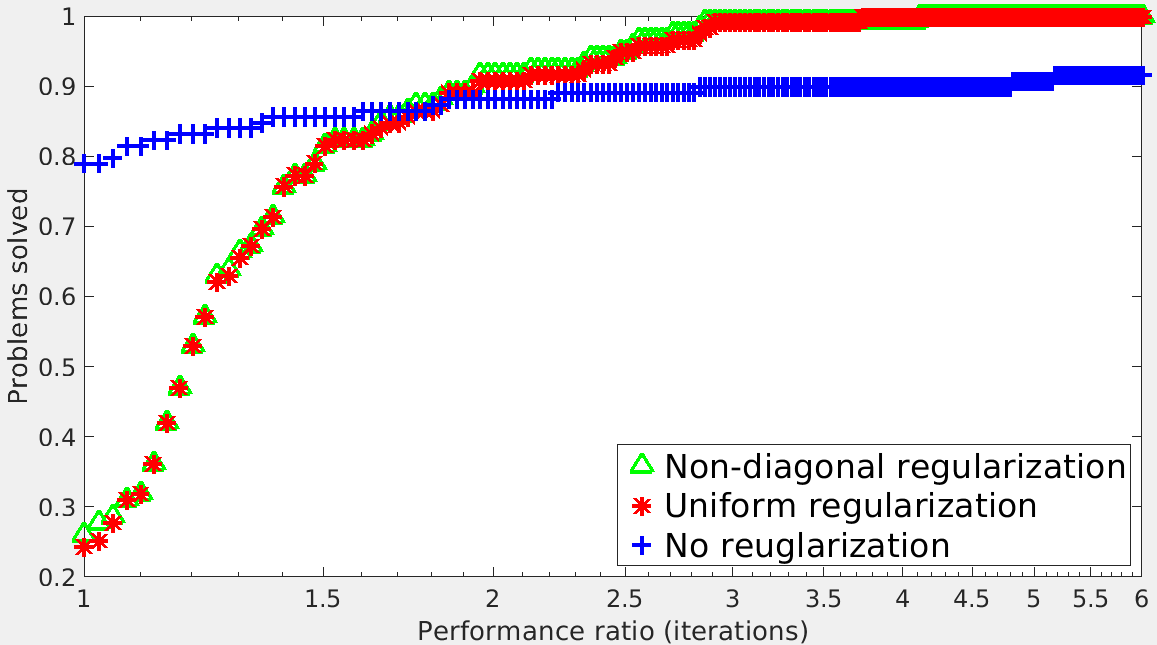}
	\caption{Performance profile with respect to iterations}
	\label{QP perf. prof. iter}
\end{subfigure}
\caption{Performance profiles over the Maros-M\'esz\'aros test set.} \label{QP perf. prof.}

\end{figure}

\noindent Following the linear programming case, we include Figure \ref{QP perf. prof.}, which contains the performance profiles, over the whole Maros-M\'esz\'aros repository of convex quadratic programming problems, of three methods; the proposed IPM with non-diagonal regularization, the IPM with uniform regularization (which was previously presented) and the same IPM but without regularization. In Figure \ref{QP perf. prof. time}, a comparison of the total time to convergence is presented, while Figure \ref{QP perf. prof. iter} contains the comparison of the total number of iterations.\\

\noindent By looking at Figure \ref{QP perf. prof.}, we can observe that as in the linear programming case, regularization seems crucial for the robustness of the method. In other words, one can observe that the IPM without regularization fails to solve 8.4\% of the problems of this test set. However, in contrast to the linear programming case, the non-regularized IPM is more efficient than the other two methods for most of the problems that it solves. This indicates that the problems in this test set are very sensitive to perturbations. The two regularization schemes, seem competitive both in terms of efficiency and robustness. In fact, the non-diagonal regularization scheme is slightly more efficient, however the difference is almost negligible. We should mention here, that the proposed tuning of the non-diagonal regularization is quite conservative. Hence, we would expect that one could improve the efficiency of such a method at the expense of its robustness. 
 
\section{Conclusions} \label{section: conlcusions}

\noindent In this paper, we derive a dynamic non-diagonal regularization scheme suitable for interior point methods. The proposed scheme is automatically tuned based on the properties of the problem, such that sufficiently large eigenvalues of the Newton system are perturbed insignificantly. The presence of non-diagonal terms in the regularization matrices allows us to introduce more sparsity in the linear system, solved to determine the Newton direction at each iteration of the interior point method. The regularization matrices can be computed expeditiously, enabling more efficient factorizations of the system matrix. The method has been implemented and the computational results demonstrate its efficiency. The results also support the claim that the proposed rule, for tuning the regularization matrices based on the properties of the problem, produces a regularization which perturbs the system almost insignificantly while maintaining numerical stability. An extension of this regularization, to interior point methods that solve the Newton system using an iterative scheme, seems natural and will be addressed in a future work. 
 \bibliography{references} 
\bibliographystyle{abbrv}
\end{document}